\documentclass[11pt]{amsart}
\usepackage{amsmath,mathtools}
\usepackage{amsmath, amsthm, amssymb, amsfonts, enumerate}
\usepackage[colorlinks=true,linkcolor=blue,urlcolor=blue]{hyperref}
\usepackage{dsfont}
\usepackage{color}
\usepackage{geometry}
\usepackage{todonotes}
\usepackage{epstopdf}
\usepackage{bbm}
\usepackage{geometry}
\usepackage[utf8]{inputenc}

\usepackage{graphicx}
\usepackage{subcaption}
\usepackage{natbib}

\usepackage{amsfonts}
\usepackage{amsfonts}
\usepackage{textcomp}
\usepackage{amssymb}
\usepackage{float}
\usepackage{tikz}
\usepackage{epsfig}
\usepackage{amsmath}
\usepackage[english]{babel}
\usepackage{a4}
\usepackage{enumerate}
\usepackage{amsmath}
\newcommand{\stkout}[1]{\ifmmode\text{\sout{\ensuremath{#1}}}\else\sout{#1}\fi}

\newtheorem{theorem}{Theorem}[section]
\newtheorem{remark}[theorem]{Remark}
\newtheorem{assumption}[theorem]{Assumption}
\newtheorem{lemma}[theorem]{Lemma}

\newtheorem{definition}[theorem]{Definition}

\def \E{\mathsf{E}}

\def \P{\mathsf{P}}
\def \R{\mathbb{R}}

\def\X{\widehat{X}}
\def\S{\widehat{S}}

\newcommand{\supp}{\text{supp}}

\definecolor{red}{rgb}{1.0,0.0,0.0}

\definecolor{blu}{rgb}{0.0,0.0,1.0}

\definecolor{gre}{rgb}{0.03,0.50,0.03}

\title[Stationary Singular Control Mean Field Games]{Stationary Discounted and Ergodic Mean Field Games of Singular Control}
\author[Cao]{Haoyang Cao}
\author[Dianetti]{Jodi Dianetti}
\author[Ferrari]{Giorgio Ferrari}

\address{H.~Cao: The Alan Turing Institute, 2QR, 96 Euston Rd, London NW1 2DB, UK}
\email{\href{hcao@turing.ac.uk}{hcao@turing.ac.uk}}
\address{J.~Dianetti: Center for Mathematical Economics (IMW), Bielefeld University, Universit\"atsstrasse 25, 33615, Bielefeld, Germany}
\email{\href{mailto:jodi.dianetti@uni-bielefeld.de}{jodi.dianetti@uni-bielefeld.de}}
\address{G.~Ferrari: Center for Mathematical Economics (IMW), Bielefeld University, Universit\"atsstrasse 25, 33615, Bielefeld, Germany}
\email{\href{mailto:giorgio.ferrari@uni-bielefeld.de}{giorgio.ferrari@uni-bielefeld.de}}

\date{\today}

\numberwithin{equation}{section}

\begin{document}

\begin{abstract} 
We study stationary mean field games with singular controls in which the representative player interacts with a long-time weighted average of the population through a discounted and an ergodic performance criterion. This class of games finds natural applications in the context of optimal productivity expansion in dynamic oligopolies. We prove existence and uniqueness of the mean field equilibria, which are completely characterized through nonlinear equations. Furthermore, we relate the mean field equilibria for the discounted and the ergodic games by showing the validity of an Abelian limit. The latter allows also to approximate Nash equilibria of -- so far unexplored -- symmetric $N$-player ergodic singular control games through the mean field equilibrium of the discounted game. Numerical examples finally illustrate in a case study the dependency of the mean field equilibria with respect to the parameters of the games.
\end{abstract}

\maketitle

\smallskip

{\textbf{Keywords}}: stationary mean field games; singular control; discounted and ergodic criterion; one-dimensional It\^o-diffusion; Abelian limit; optimal productivity expansion; $\varepsilon$-Nash equilibrium



\section{Introduction}
\label{sec:intro}

Dynamic stochastic games formalize models of competition in which different agents adjust the dynamics of underlying state variables, while optimizing a certain performance criterion that also depends on the other players' actions. The standard solution concept is that of \emph{Markov perfect equilibrium} (see \cite{Fudenberg-Tirole}), where a player's equilibrium action depends on her state and on the state variables of all the other agents at the current time. However, the complexity and difficulty of stochastic games drastically increase with the number of players. As a matter of fact -- as also argued in \cite{Adlakhaetal}, among others -- when the game faces a large number $N$ of players, the concept of Markov perfect equilibrium gives rise to two issues: computability and plausibility. On the one hand, the explicit construction of Markov perfect equilibria can be very challenging from a technical point of view, if possible at all. On the other hand, when a large number of agents play the game, the fact that each player perfectly observes and keeps track of the rivals' states can be at least debatable.

For these reasons in the past decades it has been proposed the notion of mean field equilibrium (introduced independently in \cite{HuangMalhameCaines06} and \cite{LasryLions07}; see also the two-volume book \cite{CarmonaDelarue18} for a complete overview of theory and methods) and the related one of oblivious equilibrium for infinite models, developed in \cite{Weintraubetal2008}. Both those concepts share the idea of mean field approximation from Statistical Physics: A large number of exchangeable agents playing a symmetric game with mean field interaction does not react to the actions of each single other player (as it would be in a Markov perfect equilibrium), but only to the distribution of the other agents' states, and possibly actions. 
Within this setting, stationary mean field games and stationary oblivious equilibria for infinite models (cf.\ \cite{Bardi, Cirant, Hopenhayn, Weintraubetal2011}, among others) assume that the representative player makes actions only on the basis of her own state and the long-run average state of the mass. This encompasses the following idea: In a symmetric game with a large number of players, whose state and performance criterion only depend on the distribution of opponents' state (i.e.\ an anonymous game; cf.\ \cite{JovanovicRos}), fluctuations of players' states are expected to average out, the behavior of the other agents is ``lost in the crowd'' (\cite{Gueant-etal2011}, p.\ 208), and the population's state remains roughly constant over time.  
\vspace{0.1cm}

\textbf{Our motivation and problem.} So far, most of the research on mean field games and games with oblivious equilibria has focused on the important questions of existence, uniqueness, and approximation of solutions, and analytical results in closed form have been obtained in particular models. Clearly, the difficulty of providing explicit constructions of equilibria is also due to the fact that, in general, games with mean field interaction possess an intrinsic infinite-dimensional feature, due to the dependency of the performance criterion/dynamics on the distribution of the continuum of players. It thus follows that one can hope to determine the explicit structure of mean field equilibria only in specific settings, as when the game has a linear-quadratic structure (see, e.g, \cite{Bensoussan-etal-2016} and references therein), or when the problem is stationary and the interaction is through moments of the population's state with respect to the stationary distribution (cf.\ \cite{Baseietal} and \cite{CaoGuo}, among others).

Motivated by the aim of providing the explicit construction of equilibria in games of productivity expansion, in this paper we consider a class of continuous-time stationary mean field games with singular controls. Our analysis allows to treat the mean field version of a symmetric dynamic oligopoly model, which is described as follows. In the pre-limit, each company can instantaneously increase via costly investment its productivity, which is affected by idiosyncratic noise modeling, e.g., exogenous technological shocks. In the spirit of Chapter 11 in \cite{DixitPindyck} or \cite{Bertola}, each unit of investment gives rise to a proportional cost, and investments do not need to be necessarily performed at rates; also singularly continuous actions and gulps are allowed. The operating profit function of each (exchangeable) company depends in an increasing way on the company's productivity, and it is decreasing with respect to a long-time average of a weighted mean of all other firms' asymptotic productivities. In the limit, the representative company is then expected to react to a weighted mean of the population's stationary productivity. In this paper, we actually abstract from this concrete application and study a stationary mean field game where: (i) the state variable of the representative agent is a nonnegative singularly controlled It\^o-diffusion; (ii) the interaction among players is through the reward functional, in which instantaneous profits depend on a suitable weighted average of the state process with respect to the stationary distribution; (iii) the representative agent maximizes either a discounted net profit functional or its ergodic version. Notice that a study of the ergodic mean field game is particularly relevant when one considers decisions in the context of sustainable development and management of public goods, in which it might be important to take care of the payoffs received by successive generations.
\vspace{0.1cm}

\textbf{Our contributions.} The main contributions of this work are the following. First of all, we are able to construct the unique mean field stationary equilibrium, both for a discounted and an ergodic reward functional. In both cases, the equilibrium control is of barrier-type: there exists an endogenously determined threshold $x^{\star}$ at which it is optimal to reflect the state process upward in a minimal way (i.e.\ according to a so-called Skorokhod reflection; see, e.g., Chapter 6 in \cite{Harrison}). The equilibrium stationary distribution is given by a truncated version of the speed-measure of the underlying It\^o-diffusion; that is, it coincides with the speed measure on $[x^{\star},\infty)$ and it is zero otherwise. The equilibrium triggers are characterized as the unique solutions to some nonlinear equations involving marginal profits, marginal cost of investment, and characteristic quantities of the It\^o-diffusion. Those equations can be easily solved numerically, even explicitly when, for example, the state process is a geometric Brownian motion. The approach leading to such a complete characterization of the discounted and ergodic mean field equilibria is as follows: we fix the stationary average $\theta$ of the population, and we solve one-dimensional singular stochastic control problem parametrized by $\theta$. In line with \cite{Alvarez, Alvarez2018, LZ11, JackZervos, JackJonhnsonZervos}, among others, we find that, for each given $\theta$, it is optimal to reflect the state upward at some $x^{\star}(\theta)$. We then impose that the value of $\theta$ at equilibrium is actually the one which is computed through the stationary distribution of the state process reflected at $x^{\star}(\theta)$. We prove that the resulting fixed-point problem admits a unique solution $\theta^{\star}$, which, in turn, leads to the equilibrium trigger $x^{\star}:=x^{\star}(\theta^{\star})$. It is worth observing that a byproduct of our analysis is the solution to a class of ergodic singular stochastic control problems, via exploiting a connection to optimal stopping in the spirit of \cite{Kar83}.

Second of all, we can show that the so-called \emph{Abelian limit} holds for our mean field games. This means that, when the representative agent discounts profits and costs at a rate $r$ decreasing to zero, the expected reward associated to the mean field equilibrium of the $r$-discounted problem, multiplied by $r$, converges to a suitably constant. The latter actually is the equilibrium value of the ergodic mean field game. Moreover, also the barrier triggering the equilibrium control of the discounted problem converges to that of the ergodic problem when $r\downarrow 0$. To the best of our knowledge, this is the first paper that shows the validity of the Abelian limit for a (stationary) mean field game with singular controls. The proof of such a convergence requires a careful analysis of the dependency with respect to $r$ of the equilibrium trigger and average arising in the discounted game. This is possible by analyzing the continuity with respect to $r$ of the solution to the system of equations uniquely defining the equilibria.

A natural question is whether the determined mean field equilibria approximate the corresponding symmetric $N$-player games. Moreover, in light of the Abelian limit, one can wonder whether the mean field equilibrium for the discounted problem relates to $\varepsilon$-equilibria in the ergodic symmetric $N$-player game. The study of these two questions represent the third main contribution of this work. 
We introduce ergodic and discounted $N$-player symmetric games where each player reacts to the long-time average of an increasing function of a weighted mean of the opponents' states. We then show that the mean field equilibria of the ergodic and discounted problems realize an $\varepsilon_N$-Nash equilibrium for those $N$-player games, with $\varepsilon_N$ converging to zero as $N$ goes to infinity.
Furthermore, when $N$ is large and $r$ is small, the validity of the Abelian limit allows to prove that the equilibrium of the discounted mean field game approximates a Nash equilibrium of the ergodic $N$-player singular control game. While $N$-player games with singular controls have already attracted some attention in the recent literature (see, among others, \cite{DeAFe, DianettiFerrari, FeStRi, GuoXu, Guoetal, Kwon, KwonZhang}), to our knowledge, singular control games with ergodic criterion have not yet been investigated. The previous approximation result thus sheds light on a novel class of dynamic stochastic games which naturally arise in applications.

Finally, we complement our theoretical analysis by a study of the aforementioned productivity expansion model motivating our study. Here, we assume that the productivity of the representative firm evolves as a geometric Brownian motion with negative growth rate (i.e., a positive depreciation of productivity), and that the profit function is of power type. In such a setting, the equilibrium investment triggers can be explicitly determined and numerical experiments provide the dependency of the mean field equilibria on the model's parameters. For example, consistently with the single-player optimization, one observes that, at equilibrium, larger fluctuations of the productivity leads to an increase in the ``value of waiting'', while a larger depreciation rate makes the representative company invest earlier. Further, in both the discounted and ergodic cases, the equilibrium productivity distributions are of power-type (cf.\ \cite{Gabaix}), with productivity concentrated on a right-neighborhood of the equilibrium triggers.
\vspace{0.1cm}

\textbf{Related literature.} Our work contributes to the literature on mean field games with continuous time and continuous state-space. In particular, it is placed amongst those recent works that study mean field equilibria for games with singular controls. The closest paper to ours is the recent \cite{CaoGuo}, where it is studied a stationary discounted mean field game with two-sided singular controls, and its relation to the associated $N$-player game. However, differently to our general diffusive model, in \cite{CaoGuo} the dynamics of the state process is a geometric Brownian motion and the relation to the ergodic formulation of the mean field game is not addressed. In \cite{Campietal} and \cite{GuoXu} mean field and $N$-player stochastic games for finite-fuel follower problems are studied, and the structure of equilibria is obtained. Finally, the work \cite{HorstFu} provides a careful technical analysis of the question of existence for general mean field games involving singular controls.

Closely related works are also those studying mean field games with impulsive controls and stopping times. As a matter of fact, generally speaking, singular control problems can be seen as the limit of impulsive ones when the fixed cost component vanishes (cf.\ \cite{Bensoussanetal}), and optimal timing questions can be thought of as the marginal version of singular control problems (see \cite{BalKar}). Mean field games of optimal stopping are studied in \cite{Tankov2}, \cite{Bertucci}, and \cite{Tankov1}, among others. In \cite{Tankov1}, it is proposed a relaxed-solution approach and it shown that the considered class of mean field games of stopping admits a relaxed equilibrium and that the associated value is actually unique. An application of the approach to electricity markets is then addressed in \cite{Tankov2}. In \cite{Bertucci} it is instead studied via analytical means the variational inequality associated to the mean field game of optimal timing. It is also worth noticing that mean field stationary optimal stopping problems appeared also in the economic/finance literature. Amongst many others, we refer to the early discrete-time model of firms' dynamics in \cite{Hopenhayn}, to its continuous-time version in \cite{Luttmer}, and to the competitive equilibrium model of capital structure by \cite{Miao}. In \cite{Baseietal}, it is studied nonzero-sum $N$-player and mean field stochastic games with impulse controls. Under suitable requirements, it is shown that the mean field game provides an $\varepsilon$-Nash equilibrium approximation to the $N$-player game. The results are then illustrated in a cash management problem. The work by \cite{Christensen} addresses a question similar to ours, but in the setting of impulsive mean field games. In particular, in \cite{Christensen} it is studied a class of explicitly solvable ergodic mean field control problems/games arising in optimal harvesting. However, differently to the present work, the validity of the Abelian limit and the relation to $N$-player games are not investigated there.
\vspace{0.25cm}

\textbf{Paper's structure.} The rest of the paper is organized as follows. In Section \ref{sec:setting} the probabilistic setting is introduced, while the mean field games are presented in Section \ref{sec:MFGs}. Section \ref{sec:solving} collects the results of existence and uniqueness of mean field equilibria, and Section \ref{sec:Abelianlim} derives the Abelian limit. The relation between the considered mean field games and their related symmetric $N$-player games is then discussed in Section \ref{sec:approx}. Explicit results and numerical experiments in a case study are presented in Section \ref{sec:casestudy}, while the Appendices collect the proofs of the main results of this work. 

\section{The probabilistic setting}
\label{sec:setting}

We introduce here the probabilistic setting of our study.
On a given complete filtered probability space $(\Omega,\mathcal{F},\mathbb{F}:=(\mathcal{F}_t)_t,\P)$ satisfying the usual conditions, consider a stochastic process $X^{\nu}$, whose dynamics is affected by an $\mathbb{F}$-Brownian motion $B$ and by an adapted, right-continuous, nondecreasing process $\nu$. Namely,
\begin{equation}
	\label{stateX}
	d  X^{\nu}_t = b(X^{\nu}_t) d t + \sigma(X^{\nu}_t) d B_t + d \nu_t, \quad X_{0^-} = x \in \mathbb{R}_+.
\end{equation}
Here and in the sequel, $\mathbb{R}_+ := (0, \infty)$, and $b$ and $\sigma$ are suitable drift and diffusion coefficients (see Assumption \ref{A1} below). $X^{\nu}$ models the state process of the representative agent, while $\nu$ is the control variable belonging to the set
\begin{align} 
	\mathcal{A} := & \{ \nu:\Omega \times \R_+ \to \R_+, \mbox{ $\mathbb{F}$-adapted and such that } t \mapsto \nu_t \mbox{ is a.s.} \label{setA} \\ 
	& \hspace{1cm} \mbox{nondecreasing, right-continuous and s.t. } \nu_{0^-} = 0 \}. \nonumber
\end{align}  
The nondecreasing property of the paths of $\nu$ models the fact that $\nu_t$ represents the cumulative amount of control exerted up to time $t$, such as the cumulative investment into production made up to time $t$. On the other hand, the $\mathbb{F}$-adaptedness of $\nu$ prescribes that actions should be taken on the basis of the flow of information available to the agent. 


In order to ensure that, for any given $\nu \in \mathcal{A}$ and $x\in\mathbb{R}_+$, there exists a unique strong solution to \eqref{stateX}, we make the following assumption (cf.\ Theorem 7 in Chapter V of \cite{protter2005}).
\begin{assumption}
	\label{A1}
	The coefficients $b:\mathbb{R}_+ \to \R$ and $\sigma:\mathbb{R}_+ \to \mathbb{R}_+$ are continuously differentiable. Furthermore, $b$ and $\sigma$ are (globally) Lipschitz continuous, and $\sigma\sigma'$ is locally Lipschitz.
\end{assumption}

The locally Lipschitz property of $\sigma\sigma'$, as well as the Lipschitz continuity of $b$, will be needed in our subsequent analysis (cf.\ \eqref{statehatX} and Proof of Lemma \ref{lemma phi continuous}, respectively). We denote the solution to \eqref{stateX} by $X^{x,\nu}$, and, when needed, we use the notation $\E[g(X^{x,\nu})] = \E_x[g(X^{\nu})]$, for any Borel-measurable and integrable function $g$. 

Furthermore, the uncontrolled diffusion process $X:=X^0$ solving                               
\begin{equation}
	\label{stateXbar}
	d {X}_t = b(X_t) d t + \sigma(X_t) d B_t, \quad X_0 = x \in \mathbb{R}_+,
\end{equation}
is nondegenerate, and for any $x_o \in \mathbb{R}_+$ there exists $\varepsilon > 0$ (depending on $x_o$) such that 
\begin{equation}
	\label{LIX}
	\int_{x_o - \varepsilon}^{x_o + \varepsilon} \frac{1 + |b(z)|}{\sigma^2(z)}\ dz < + \infty.
\end{equation}
The latter guarantees that $X$ is a regular diffusion. That is, starting from $x \in \mathbb{R}_+$, $X$ reaches any other $y \in \mathbb{R}_+$ in finite time with positive probability. Finally, to stress the dependency of $X$ on its initial value, from now on we may write $X^{x}$. 

In our subsequent analysis, an important role will be also played by the one-dimensional It\^o-diffusion $\widehat X$ evolving as
\begin{equation}
	\label{statehatX} 
	d \widehat X_t = \big[b(\widehat X_t) + (\sigma \sigma')(\widehat X_t)\big] d t + \sigma(\widehat X_t) d \widehat B_t, \quad \widehat X_0 = x \in \mathbb{R}_+,
\end{equation}
for some one-dimensional $\mathbb{F}$-Brownian motion $\widehat B$.

Notice that, under Assumption \ref{A1}, there exists a unique strong solution to \eqref{statehatX}, up to a possible explosion time. Moreover, one has that for any $x_o \in \mathbb{R}_+$ there exists $\varepsilon > 0$ such that 
\begin{equation}
	\label{LIhatX}
	\int_{x_o - \varepsilon}^{x_o + \varepsilon} \frac{1 + |b(z)| + |\sigma \sigma'(z)|}{\sigma^2(z)}d  z < + \infty,
\end{equation}
ensuring that $\widehat X$ is a regular diffusion as well. In order to highlight the dependence of $\X$ on its initial value, in the following we write $\X^x$, when needed, and we denote by $\widehat{\E}_x$ the expectation under the measure $\widehat{\P}_x$ on $(\Omega, \mathcal{F})$ such that $\widehat{\P}_x[\,\cdot\,]:=\P[\,\cdot\,|\,\widehat X_0=x]$.


\subsection{On the diffusions $X$ and $\widehat{X}$: Characteristics and requirements}
\label{sec:diffusions}   

In this subsection we recall useful basic characteristics of the diffusion processes $X$ and $\X$. We refer to Chapter II in \cite{Borodin Salminen} for further details.

The infinitesimal generator of the uncontrolled diffusion $X$ is denoted by $\mathcal{L}_{X}$ and is defined as the second-order differential operator
\begin{equation} 
	\label{eq LX} 
	(\mathcal{L}_{X} f) \,(x):=\frac{1}{2}\sigma^2(x)f''(x)+b(x)f'(x), \quad f \in C^2(\mathbb{R}_+), \quad x \in \mathbb{R}_+. 
\end{equation} 
On the other hand, the infinitesimal generator $\X$ is denoted by $\mathcal{L}_{\X}$ and is such that
\begin{align}
	\label{eq LXhat}
	(\mathcal{L}_{\X} f)\,(x):=\frac{1}{2}\sigma^2(x)f''(x)+(b(x)+\sigma(x)\sigma'(x))f'(x),\quad f\in C^2({\mathbb{R}_+}), \quad x\in\mathbb{R}_+.
\end{align}

For $r > 0$, we introduce $\psi_r$ and $\phi_r$ as the fundamental solutions to the ordinary differential equation (ODE),
\begin{align}
	\label{ODE}
	\mathcal{L}_X u(x)-ru(x)=0,\qquad x\in\mathbb{R}_+,
\end{align}
and we recall that they are strictly increasing and decreasing, respectively. For an arbitrary $x_o \in \mathbb{R}_+$ we also denote by 
$$S'(x) := \exp\left(- \int_{x_o}^x \frac{2b(z)}{\sigma^2(z)}\ d  z \right), \qquad x\in\mathbb{R}_+, $$
the derivative of the scale function of $X$, and we observe that
the derivative of the speed measure of $X$ is given by $m'(x):=\frac{2}{\sigma^2(x)\, S'(x)}$.
Together with the killing measure, scale function and speed measure represent the basic characteristics of any diffusion process. In particular, $S$ is related to the drift of the diffusion and, more specifically, to the probability of the diffusion leaving an interval either from its left or right endpoint. On the other hand, it can be shown that the transition probability of a regular diffusion is absolutely continuous with respect to the speed measure. 

Throughout this paper we assume that 
\begin{assumption}
	\label{ass:X-rec}
	$$\int_{a}^{\infty} m'(y) d  y < \infty, \quad \text{for any}\quad a > 0.$$
\end{assumption}
This requirement guarantees that the process $X$ reflected upward at a level $a>0$ is positively recurrent. 

Moreover, when $r-b'(x)\geq r_o>0$ for $x\in {\mathbb{R}_+}$, any solution to the ODE
\begin{align}
	\label{ODE2}
	\mathcal{L}_{\X} u(x)-(r-b'(x))u(x)=0,\qquad x\in\mathbb{R}_+,
\end{align}
can be written as a linear combination of the fundamental solutions $\widehat{\psi}_r$ and $\widehat{\phi}_r$, which are strictly increasing and decreasing, respectively. 
Finally, letting $x_o \in \mathbb{R}_+$ to be arbitrary, we denote by 
$$ \S'(x) := \exp\left(- \int_{x_o}^x \frac{2b(z) + 2 \sigma (z) \sigma'(z)}{\sigma^2(z)}\ d  z \right), \qquad x\in\mathbb{R}_+, $$ 
the derivative of the scale function of $\X$, and by $\widehat{m}'(x) := \frac{2}{\sigma^2(x)\, \S'(x)}$ the density of its speed measure.
One can easily check that the scale functions and speed measures of $X$ and $\X$ are related through $\S'(x)=S'(x)/\sigma^2(x)$ and $\widehat{m}'(x)=2/S'(x)$, for $x \in \mathbb{R}_+$.

Concerning the boundary behavior of the real-valued It\^o-diffusions $X$ and $\widehat X$, in the rest of this paper \textbf{we assume that $0$ and $+\infty$ are natural boundaries} for those two processes. In particular, this means that $0$ and $\infty$ are unattainable in finite time and that, for each $r > 0$, we have
\begin{equation}
	\label{psiphiproperties1}
	\lim_{x \downarrow 0}\psi_r(x) = 0,\,\,\,\,\lim_{x \downarrow 0}\phi_r(x) = + \infty,\,\,\,\,\lim_{x \uparrow \infty}\psi_r(x) = + \infty,\,\,\,\,\lim_{x \uparrow \infty}\phi_r(x) = 0,
\end{equation}
\begin{equation}
	\label{psiphiproperties2}
	\lim_{x \downarrow 0}\frac{\psi_r'(x)}{S'(x)} = 0,\,\,\,\,\lim_{x \downarrow 0}\frac{\phi_r'(x)}{S'(x)} = -\infty,\,\,\,\,\lim_{x \uparrow \infty}\frac{\psi_r'(x)}{S'(x)} = + \infty,\,\,\,\,\lim_{x \uparrow \infty}\frac{\phi_r'(x)}{S'(x)} = 0.
\end{equation}
Also, when $r-b'(x) \geq r_o >0$ for each $x \in \mathbb{R}_+$, we have 
\begin{equation}
	\label{psiphiproperties1bis}
	\lim_{x \downarrow 0}\widehat{\psi}_r(x) = 0,\,\,\,\,\lim_{x \downarrow 0}\widehat{\phi}_r(x) = + \infty,\,\,\,\,\lim_{x \uparrow \infty}\widehat{\psi}_r(x) = + \infty,\,\,\,\,\lim_{x \uparrow \infty}\widehat{\phi}_r(x) = 0,
\end{equation}
\begin{equation}
	\label{psiphiproperties2bis}
	\lim_{x \downarrow 0}\frac{\widehat{\psi}_r'(x)}{\S'(x)} = 0,\,\,\,\,\lim_{x \downarrow 0}\frac{\widehat{\phi}_r'(x)}{\S'(x)} = -\infty,\,\,\,\,\lim_{x \uparrow \infty}\frac{\widehat{\psi}_r'(x)}{\S'(x)} = + \infty,\,\,\,\,\lim_{x \uparrow \infty}\frac{\widehat{\phi}_r'(x)}{\S'(x)} = 0.
\end{equation}

Furthermore, \textbf{we require that} 
$$\lim_{x \downarrow 0}\phi_r'(x)=-\infty \quad \text{and} \quad \lim_{x \uparrow \infty}\psi_r'(x)=\infty.$$ 
Then, by arguing as in the second part of the proof of Lemma 4.3 in \cite{AlvarezMat}, one can show that, under our conditions on $X$ and $\X$, one has $\widehat{\phi}_r=-\phi_r'$ and  $\widehat{\psi}_r=\psi_r'$.

Finally, the following useful equations hold for any $0<a<b<\infty$:
\begin{equation}
	\label{psiphiproperties3}
	\left\{ \begin{array}{ll}
		\displaystyle \frac{\widehat{\psi}_r'(b)}{\S'(b)} - \frac{\widehat{\psi}_r'(a)}{\S'(a)}= \int_{a}^{b}\widehat{\psi}_r(y)(r - b'(y))\widehat{m}'(y) d  y, \vspace{0.25cm} \\
		\displaystyle \frac{\widehat{\phi}_r'(b)}{\S'(b)} - \frac{\widehat{\phi}_r'(a)}{\S'(a)}= \int_{a}^{b}\widehat{\phi}_r(y)(r - b'(y))\widehat{m}'(y) d  y.	
	\end{array} \right.
\end{equation} 

We conclude this discussion by noticing that all the requirements on $X$ (and, consequently, on $\X$) assumed so far are satisfied, for example, by the relevant cases in which $X$ is a geometric Brownian motion with drift $b(x)=-\delta x$, $\delta>0$, or an affine mean-reverting dynamics with drift $b(x)=\kappa(\lambda - x)$ and volatility $\sigma(x)=\sigma x$, for positive $\kappa, \lambda, \sigma$.


\section{The stationary mean field games}
\label{sec:MFGs}

In this section we introduce the stationary mean field games (MFGs) that will be the object of our study. 
For any $\nu \in \mathcal{A}$, $x \in \mathbb{R}_+$, and $\theta \in \mathbb{R}_+$ we consider the \emph{discounted} expected profit
\begin{equation}
	\label{profit-d}
	J (x,\nu, \theta; r) := \E_x\bigg[ \int_0^{\infty} e^{- r s} \pi(X^{\nu}_s, \theta) ds -  \int_0^{\infty} e^{- r s} d \nu_s \bigg], \quad r>0,
\end{equation} 
as well as the \emph{ergodic} expected profit
\begin{align}
	\label{profit-e} 
	& G (x, \nu, \theta) :=   
	\displaystyle \limsup_{T\uparrow \infty}\frac{1}{T}\E_x\bigg[ \int_0^{T} \pi(X^{\nu}_s, \theta) ds -  \nu_T \bigg]. 
\end{align}
In \eqref{profit-d} and \eqref{profit-e}, $\pi:\R^2 \to [0,\infty)$ is an instantaneous profit function and the control processes are picked from the admissible classes
\begin{align}\label{def admiss controls} 
	& \mathcal{A}_d := \Big\{ \nu \in \mathcal{A}:\, \E\Big[\int_0^{\infty} e^{- r s} d \nu_s\Big] < \infty \Big\} \quad \text{and} \quad \mathcal{A}_{e} :=  \Big\{ \nu \in \mathcal{A}:\, \E\big[\nu_T\big] < \infty \,\, \forall T >0\Big\},
\end{align}
respectively. In order to simplify notation, in the sequel we shall omit the dependency on $r$ of the set $\mathcal{A}_d$, which in fact will be clear from the context. Furthermore, in \eqref{profit-d} the integral $\int_0^{t}(\,\cdot\,) d\nu_s$ is intended in the Lebesgue-Stieltjes sense as $\int_{[0,t]}(\,\cdot\,) d\nu_s$, thus including a possible initial jump of $\nu$ of amplitude $\nu_0$.


Next, for a probability measure $\mu$ on $\R$, i.e.\ $\mu \in \mathcal{P}(\R_+)$, we define 
\begin{equation}
	\label{def:theta}
	\theta(\mu) :=  F\bigg( \int_{\mathbb{R}_+} f(x) \mu(d  x) \bigg) \in  [0, \infty], 
\end{equation}
where the functions $F$ and $f$ are a strictly increasing nonnegative functions. 

In the mean field games defined through the next Definitions \ref{def equilibrium discounted}  and \ref{def equilibrium ergodic}, the term $\theta=\theta(\mu)$ appearing in \eqref{profit-d} and \eqref{profit-e} describes a suitable mean with respect to the stationary distribution $\mu=\P_{X^{\nu}_{\infty}}$ of the optimally controlled state process $X^{\nu}$ (provided that one exists). For example, if $X^{\nu}$ describes the productivity of the representative company, then $\mu$ provides the distribution of the asymptotic productivity, and its weighted average -- with weight function $f$ -- defines a price index through the function $F$ (cf.\ Remark \ref{rem:ISO} below).

In the sequel, we focus on the following definition of MFG equilibria.  

\begin{definition}[Equilibrium of the discounted MFG]
	\label{def equilibrium discounted} 
	For $r>0$ and $x\in \mathbb{R}_+$, a couple $(\nu^{\star,r}, \theta_r^\star) \in \mathcal{A}_d \times \mathbb{R}_+$ is said to be an equilibrium of the discounted MFG for the initial condition $x$ if 
	\begin{enumerate}
		\item $J (x,\nu^{\star,r}, \theta_r^\star ; r) \geq J (x,\nu, \theta_r^\star ; r) $, for any $\nu \in \mathcal{A}_d$;
		\item The optimally controlled process $X^{x,\star}:=X^{x,\nu^{\star,r}}$ admits a limiting distribution  $\P_{X^{x,\star}_\infty}$ satisfying $\theta_r^\star = \theta (\P_{X^{x,\star}_\infty})$.
	\end{enumerate}
\end{definition} 

\begin{definition}[Equilibrium of the ergodic MFG]
	\label{def equilibrium ergodic} 
	For $x\in \mathbb{R}_+$, a couple $(\nu^{\star, e}, \theta_e^\star) \in \mathcal{A}_e \times \mathbb{R}_+$ is said to be an equilibrium of the ergodic MFG for the initial condition $x$ if 
	\begin{enumerate}
		\item $G (x,\nu^{\star,e}, \theta_e^\star ) \geq G (x,\nu, \theta_e^\star ) $, for any $\nu \in \mathcal{A}_e$;
		\item The optimally controlled process $X^{x,\star}:=X^{x,\nu^{\star,e}}$ admits a limiting distribution $\P_{X^{x,\star}_\infty}$ satisfying $\theta_e^\star = \theta (\P_{X^{x,\star}_\infty})$.
	\end{enumerate}
\end{definition}


We enforce the following structural conditions on the running profit and weight function.
\begin{assumption} 
	\label{A2} The running profit $\pi:\mathbb{R}_+^2 \to [0,+ \infty)$ belongs to $C^2(\mathbb{R}_+^2)$. Furthermore, 
	\begin{itemize}
		\item[($\pi$-i)] $\pi(\cdot,\theta)$ is concave and nondecreasing for any $\theta \in \mathbb{R}_+$;
		\item[($\pi$-ii)] $\pi$ has strictly decreasing differences; that is, $\pi_{x\theta} (x,\theta) < 0$ for any $(x,\theta) \in \mathbb{R}_+^2$;
		\item[($\pi$-iii)] for any $x\in \mathbb{R}_+$,
		$$\lim_{\theta \downarrow 0} \pi_x(x,\theta)=+\infty \quad \text{and} \quad \lim_{\theta \uparrow \infty} \pi_x(x,\theta)=0;$$
	 \item[($\pi$-iv)] for any $0<a<b< \infty$ there exists a function $h^{a,b}:\mathbb{R}_+ \to \mathbb{R}_+$ such that
	$$ 
	|\pi_{x\theta}(x,\theta)|\,  \widehat{m}'(x) \leq h^{a,b}(x), \ \text{ for any } x \in \mathbb{R}_+, \ \theta \in (a,b), \quad \text{and} \quad h^{a,b} \in \mathbb{L}^1(\kappa,\infty)  \text{ for any } \kappa \in \mathbb{R}_+.
	$$
	\end{itemize}
	
	Moreover, the weight functions $F$ and $f$ appearing in \eqref{def:theta} satisfies:
	\begin{itemize}
		\item[(i)] The functions $F,f:[0,\infty) \to [0,\infty)$ are continuously differentiable with $F', f' >0$ and, for $\beta \in (0,1)$ and a constant $C>0$, they satisfy the  growth conditions: 
		\begin{align*}
			f(x) &\leq C (1+|x|^\beta), \\
			F(x) & \leq C(1+|x|^{\frac{1}{\beta}} ), \\
			|F(y) - F(x)| & \leq C(1+|x|+|y|)^{\frac{1}{\beta} -1 }|y-x|, 
		\end{align*}
		for any $y, x \in \mathbb R_+$;
		\item[(ii)] $\lim_{y\uparrow \infty} F(y) = + \infty$ and $\lim_{y\uparrow \infty} f(y) = + \infty$.
	\end{itemize} 
\end{assumption}  

\begin{remark}
\label{rem:ISO}
With regard to a stationary mean field formulation of a game of productivity expansion, a benchmark example of running profit function and average satisfying Assumption \ref{A2} are
$$\pi(x,\theta):=x^{\beta}\,\theta^{-(1 + \beta)},\quad \theta:=\theta(\mu)=\bigg( \int_{\mathbb{R}_+}x^{\beta}\mu(d x)\bigg)^\frac{1}{\beta},  \quad \beta \in (0,1).$$  
Such a form of interaction can be obtained from the so-called isoelastic demand obtained from Spence-Dixit-Stiglitz preferences (see, e.g., footnote 5 in \cite{Achdouetal} for such a derivation).
\end{remark}


\section{Existence, uniqueness, and characterization of the mean field equilibria} 
\label{sec:solving}

In this section, the discounted MFG problem and the ergodic MFG problem are solved. In particular, existence and uniqueness of equilibria is shown by charaterizing the equilibria in terms of the unique solution to systems of nonlinear equations.

\subsection{On the discounted stationary MFG}
\label{sec:discMFG}

In order to deal with the discounted MFG problem for a fixed discount factor $r>0$, we make the following additional requirement (see also \cite{JackJonhnsonZervos}, \cite{Kwon}, among others). 
\begin{assumption}
	\label{assumption discounted}\
	\begin{itemize} 
		\item[(i)] For each $x\in \mathbb{R}_+$ we have $r-b'(x)\geq 2c>0$, for a constant $c>0$;
		\item[(ii)] For any $\theta \in \R_+$, there exists $\widehat{x}_r(\theta) \in \R_+$ such that
		$$ \pi_x(x,\theta) - r + b'(x)\left\{ \begin{array}{ll}
		< 0, \quad	& x > \widehat{x}_r(\theta), \\
		= 0, \quad 	& x = \widehat{x}_r(\theta), \\
		> 0, \quad 	& x < \widehat{x}_r(\theta).
		\end{array} \right. $$
	\end{itemize}
\end{assumption}
Condition ($i$) above guarantees that the discount rate is (uniformly) larger than the marginal growth rate of the diffusion $X$. It is automatically satisfied in the particular cases in which $X$ is a geometric Brownian motion with drift $b(x)=-\delta x$, $\delta>0$, or it is an affine mean-reverting process with drift $b(x)=\kappa(\lambda-x)$, $\kappa,\lambda>0$ (and volatility $\sigma(x)=\sigma x$, $\sigma>0$). Moreover, bearing in mind the mean field game of productivity expansion discussed in the introduction, Condition ($ii$) in Assumption \ref{assumption discounted} ensures the following: The marginal running profit $\pi_x$, net of the ``user cost of capital'' $r-b'$, changes sign at most once. Such a requirement guarantees that the mean field equilibrium is of threshold type, as in fact, for any given $\theta$, it should not be profitable to increase productivity via costly investment when $\pi_x(x,\theta) - r + b'(x) <0$. This is formalized in the following theorem, whose proof can be found in Appendix \ref{proof:Thm-discMFG}.

\begin{theorem}
	\label{thm discMFG-ex}
	Let $r>0$, and let Assumptions \ref{A1}, \ref{ass:X-rec}, \ref{A2},  and \ref{assumption discounted} hold. For any $x \in \mathbb{R}_+$, there exists a unique equilibrium  $(\nu^{\star,r}, \theta_r^{\star})$ of the discounted MFG. Moreover, $\nu^{\star ,r}$ makes the state process reflected upward at the barrier $x_r^{\star} < \widehat{x}_r(\theta_r^{\star})$, and the couple  $(x_r^{\star}, \theta_r^{\star})$ is determined as the unique solution to the system 
	\begin{equation}\label{eq system discounted}
		\int_{x_r^{\star}}^{\infty} \widehat{\phi}_r (y) \big(\pi_x(y,\theta_r^{\star})- r + b'(y)\big) \widehat m'(y) d  y =0 \quad \text{and} \quad \int_{x_r^{\star}}^{\infty} (f(y) - F^{-1}(\theta_r^{\star})) m'(y) d  y = 0.
	\end{equation} 
\end{theorem}


\subsection{On the ergodic stationary MFG} 
\label{sec:erg}

Our analysis of the ergodic MFG problem is subject to the following requirements, which are consistent to those in Assumption \ref{assumption discounted} when $r=0$.

\begin{assumption}
	\label{assumption ergodic} \ 
	\begin{itemize} 
		\item[(i)] For each $x\in \mathbb{R}_+$ we have $b'(x) < -2c < 0$, for a constant $c>0$; 
		\item[(ii)] For any $\theta \in \R_+$, there exists  $\widehat{x}_0(\theta) \in \R_+$ such that
		$$ \pi_x(x,\theta) + b'(x) \left\{ \begin{array}{ll}
		< 0,	& x > \widehat{x}_0(\theta), \\
		= 0,	& x = \widehat{x}_0(\theta), \\
		> 0,	& x < \widehat{x}_0(\theta).
		\end{array} \right. $$
	\end{itemize}
\end{assumption}
Notice that the condition $b'(x) < -2c < 0$ is easily seen to be verified in the relevant cases of $X$ being a geometric Brownian motion and a mean-reverting affine process with drift $b(x)=\kappa(\lambda-x)$, $\kappa,\lambda>0$, (and volatility $\sigma(x)=\sigma x$, $\sigma>0$).

Recall now the mean field game problem with ergodic net profit given by \eqref{profit-e}, together with its notion of solution given in Definition \ref{def equilibrium ergodic}. For each $\theta > 0$, set
\begin{equation}
	\label{eq control proble ergodic} 
	\lambda^{\star}(\theta):= \sup_{\nu \in \mathcal{A}_e } \displaystyle \limsup_{T\uparrow \infty}\frac{1}{T}\E_x\bigg[ \int_0^{T} \pi(X^{\nu}_s, \theta) d  s -  \nu_T \bigg] 
\end{equation}

The next result provides a complete characterization of the ergodic mean field equilibrium. Its proof can be found in Appendix \ref{sec:proof-Thm-ergMFG}.

\begin{theorem}
	\label{thm ergMFG-ex} 
	Let Assumptions \ref{A1}, \ref{ass:X-rec}, \ref{A2}, and \ref{assumption ergodic} hold. For any $x \in \mathbb{R}_+$, there exists a unique equilibrium $(\nu^{\star, e}, \theta_e^{\star})$ of the ergodic MFG. 
	Moreover, the process $\nu^{\star,e}$ reflects the state process at the barrier ${x}_e^\star < \widehat{x}_0({\theta}_e^\star)$, and the
	couple $({x}_e^\star, {\theta}_e^\star)$ is determined as the unique solution to the system   
	\begin{equation}\label{eq system for the ergodic MFG}
		\int_{x_e^\star}^{\infty} \widehat{\phi}_0 (y) \big(\pi_x(y,\theta_e^\star) + b'(y)\big) \widehat m'(y) d  y =0 \quad \text{and} \quad \int_{x_e^\star}^{\infty} (f(y) - F^{-1}(\theta_e^\star)) m'(y) d  y = 0.
	\end{equation}
	Finally, the value of the ergodic MFG  at equilibrium is given by   
	\begin{equation}\label{eq value MFG ergodic cost}         
		\lambda^\star (\theta_e^\star)  = b(x_e^\star) +  \pi (x_e^\star, \theta_e^\star).
	\end{equation}
\end{theorem}

%

\section{Connecting discounted and ergodic MFGs: The Abelian limit}
\label{sec:Abelianlim} 

A natural question is whether the mean field equilibrium and the relative equilibrium value of the discounted game can be related to those of the ergodic game in the limit $r\downarrow 0$. In this section we provide a positive answer to the previous question by showing the validity of the so-called Abelian limit for the equilibrium value of the discounted game. Moreover, we also prove convergence of the equilibrium boundary of the discounted game towards that of the ergodic game. Although similar results are known in the literature on singular stochastic control problems (cf.\ \cite{AlvarezHening}, \cite{Kar83}, \cite{Weerasinghe2007}), to our knowledge they appear here for the first time within this literature in the mean field context.

The main idea of the subsequent analysis is to show suitable regularity, with respect to the discount factor $r$ in a neighborhood of $0$, of the solutions to the systems of equations provided in Theorems \ref{thm discMFG-ex} and \ref{thm ergMFG-ex}, which in fact completely characterize the MFG equilibria. Throughout this section, \textbf{we let Assumptions \ref{A1}, \ref{ass:X-rec}, \ref{A2}, and \ref{assumption ergodic} hold}, so that Assumption \ref{assumption discounted} is also satisfied for any $r>0$.

Let then $c>0$ be as in Assumption \ref{assumption ergodic}, and define the functions 
\begin{align*}
	\Pi(y,\theta;r)&:=\widehat m'(y) \big(\pi_x(y,\theta)- (r-b'(y))\big) 
	&x,\, \theta >0, \ r\in(-c,1),\\
	K(x,\theta;r)&:=\int_{x}^{\infty} \widehat{\phi}_r (y) \Pi(y,\theta,r) d  y, &x,\, \theta >0, \ r\in (-c,1), \\
	\widehat{K}(x,\theta;r)&:=K(x,\theta;r)/\widehat{\phi}_r (x), 
	& x,\,\theta >0, \ r\in (-c,1),\\
	\widehat{G}(x,\theta)&:=\int_{x}^{\infty} (f(y) - F^{-1}(\theta)) m'(y) d  y 
	& x,\, \theta >0. 
\end{align*}
Define next $\Phi:\R_+^2  \times (-c,1) \to \R^2$ by setting
\begin{equation}
	\label{eq:PHI}
	\Phi(x,\theta;r):=(\widehat{K}(x,\theta;r),\widehat{G}(x,\theta)). 
\end{equation}
The function $\Phi$ describes the system of equations determining the MFG equilibria of the discounted problem and of the ergodic problem.
Indeed, for each $r> 0$, since $\widehat{\phi}_r>0$, we have $K(x,\theta;r)=0$ if and only if $\widehat{K}(x,\theta;r)=0$; hence, according to Theorem \ref{thm discMFG-ex}, for each $r > 0$, there exists a unique $(x_r^{\star},\theta_r^\star)$ such that $\Phi(x_r^{\star},\theta_r^\star;r)=0$. Analogously, according to Theorem \ref{thm ergMFG-ex}, there exists a unique $(x_e^{\star},\theta_e^\star)$ such that $\Phi(x_e^{\star},\theta_e^\star;0)=0$. 

Clearly, continuity of $\Phi$ is a necessary ingredient for the previously discussed convergence of the equilibrium of the discounted MFG towards that of the ergodic MFG. This is accomplished in the next technical lemma, whose proof is postponed to Appendix \ref{proof:lemma phi continuous}.

\begin{lemma}
	\label{lemma phi continuous}
	The function $\Phi:\R_+^2 \times (-c,1) \to \R^2$ is continuous.
\end{lemma}

For each $r\in (0,1]$, denote now by $V(x,r;\theta_r^\star)$ the equilibrium value of the MFG with discount factor $r$. We are then in the condition of stating (and proving in Appendix \ref{proof:thm: abelian}) the main result of this section.

\begin{theorem}
	\label{thm: abelian}
	For any $x \in \mathbb{R}_+$, one has
	$$
	\lim_{r\downarrow 0}(x^{\star}_r, \theta^{\star}_r) = (x^{\star}_e, \theta^{\star}_e)\quad \text{and}\quad \lim_{r\downarrow 0} rV(x,\theta^{\star}_r;r)= \lambda^\star(\theta^{\star}_e).$$
\end{theorem}


\section{Mean-field vs.\ $N$-player: approximation results} 
\label{sec:approx} 

In the previous sections, we have established existence and uniqueness of the solutions to both the discounted and ergodic MFGs with singular controls. 
Here we provide the connection of these mean field solutions to symmetric $N$-player games. In particular, we show that each mean field solution approximates the Nash equilibrium of a suitable $N$-player game. Furthermore, by exploiting the Abelian limit, we find that the mean field equilibrium of the discounted game realizes an $\varepsilon$-Nash equilibrium for the $N$-player ergodic game, when $N$ is large and $r$ is small. These results have the two following implications: on the one hand, they shed light on the ``closeness'' of $N$-player discounted games with the -- so far unexplored -- $N$-player ergodic games, when $N$ is large and $r$ is small; on the other hand, they provide an operative way of constructing approximate equilibria.

Throughout this section, \textbf{we let Assumptions \ref{A1}, \ref{ass:X-rec}, \ref{A2}, and \ref{assumption ergodic} hold}, so that Assumption \ref{assumption discounted} is also satisfied for any $r>0$. 

The $N$-player games are described as follows. Let the filtered probability space $(\Omega, \mathcal F,\mathbb F=(\mathcal F_t)_{t\geq0},\P)$ support a standard Brownian motion $W$, and a sequence $(W^i)_{i\in \mathbb{N}}$ of independent $\mathbb{F}$-Brownian motions, independent from $W$.  
Suppose also that the filtered probability space is rich enough to allow for a sequence $(\xi^i)_{i \in \mathbb{N}}$ of i.i.d.\ square-integrable $\mathbb{R}_+$-valued  $\mathcal{F}_0$-random variables, independent from $W$ and $(W^i)_{i \in \mathbb{N}}$, and with distribution $\mu_0$.
For each $i\in\{1,\dots, N\}$, player $i$ chooses an (open-loop) strategy $\nu^i \in \mathcal{A}$ in order to control its state process $X^{i,\nu^i}$, which evolves according to 
\begin{equation}
	\label{eq n-dyn-i}
	dX^{i,\nu^i}_{t}=b(X^{i,\nu^i}_{t})dt+\sigma(X^{i,\nu^i}_{t})dW^i_t+d\nu^i_t, \quad  X_{0^-}^{i,\nu^i}=\xi^i.
\end{equation}

For strategies $\nu^i \in \mathcal{A}$, we denote by $\boldsymbol{\nu}^{-i}=(\nu^1,\dots,\nu^{i-1},\nu^{i+1},\dots,\nu^N)$ the vector of strategies picked by player $i$'s opponents, and we define profile strategies by $(\nu^i,\boldsymbol{\nu}^{-i}) := (\nu^1,\dots,\nu^N)$.
We set 
\begin{equation}\label{eq definition of theta N}
	\theta^{N}_{\boldsymbol{\nu}^{-i}}:= \lim_{t \to \infty} \frac{1}{t} \int_0^t F\Big( \frac{1}{N-1} \sum_{j \ne i} f(X_s^{j, \nu^j}) \Big) ds,
\end{equation}
and, for $q\in\{d,e\}$ and $\mathcal{A}_q$ as in \eqref{def admiss controls}, introduce the sets:
$$
\widehat{\mathcal{A}}^{N-1}_q:=\big\{\boldsymbol{\nu}^{-i} \in \mathcal{A}_q^{N-1}:\,\theta^{N}_{\boldsymbol{\nu}^{-i}} \text{ exists finite a.s.}\big\}.
$$ 
Then, the ergodic and the discounted payoffs of player $i$ arising from playing $\nu^i\in \mathcal{A}_e$ and $\nu^i\in \mathcal{A}_d$, respectively, are given by
\begin{equation} 
	\label{eq n-payoff-i ergodic}
	G^i(\nu^i,\boldsymbol{\nu}^{-i}):=\limsup_{T\to\infty}\frac{1}{T}\E\biggl[\int_0^T\pi\left(X^{i,\nu^i}_{t},\theta^{N}_{\boldsymbol{\nu}^{-i}}\right) dt-\nu^i_T \biggl], \quad \boldsymbol{\nu}^{-i} \in \widehat{\mathcal{A}}^{N-1}_e,
\end{equation} 
and  
\begin{equation} 
	\label{eq n-payoff-i discounted}
	J^i(\nu^i,\boldsymbol{\nu}^{-i};r):=\E\biggl[\int_0^\infty e^{-rt} \pi\left(X^{i,\nu^i}_{t},\theta^{N}_{\boldsymbol{\nu}^{-i}}\right) dt- \int_0^\infty e^{-rt}d\nu^i_t\biggl], \quad \boldsymbol{\nu}^{-i} \in \widehat{\mathcal{A}}^{N-1}_d.
\end{equation}

\begin{definition}[$\varepsilon$-Nash Equilibrium]
\label{def:epsilonNash}
	For $\varepsilon>0$, 
	\begin{enumerate}
		\item $\boldsymbol{\nu}^{\star}=(\nu^{{\star},1},\dots,\nu^{{\star},N}) \in {\mathcal{A}}_e^N$ is called $\varepsilon$-Nash equilibrium ($\varepsilon$-NE) of the ergodic $N$-player game if for any $i=1,\dots, N$ we have $\boldsymbol{\nu}^{ \star,-i} \in \widehat{\mathcal{A}}_e^{N-1}$ and
		\[ 
		G^i(\nu^{\star,i},\boldsymbol{\nu}^{\star,-i})\geq G^i(\nu^{i};\boldsymbol{\nu}^{\star,-i})-\varepsilon,\quad \nu^i\in \mathcal{A}_e;
		\] 
		\item  $\boldsymbol{\nu}^{\star}=(\nu^{{\star},1},\dots,\nu^{{\star},N}) \in \mathcal{A}_d^N$ is called $\varepsilon$-Nash equilibrium ($\varepsilon$-NE) of the discounted $N$-player game if for any $i=1,\dots, N$ we have $\boldsymbol{\nu}^{ \star,-i} \in \widehat{\mathcal{A}}_d^{N-1}$ and
		\[ 
		J^i(\nu^{\star,i},\boldsymbol{\nu}^{\star,-i};r)\geq J^i(\nu^{i};\boldsymbol{\nu}^{\star,-i};r)-\varepsilon,\quad \nu^i\in \mathcal{A}_d. 
		\] 
	\end{enumerate}
\end{definition}

In order to approximate Nash equilibria when the initial conditions for the SDEs \eqref{eq n-dyn-i} are random variables, for any $\theta >0 $ we define the profit functionals for the mean field game problems:
\begin{align}\label{def profit functional MFG random initial condition}
	G(\nu, \theta)&:= \int_{\R_+} G(x,\nu,\theta) \mu_0(dx), \quad \nu \in \mathcal{A}_e,
	\\ \notag
	J(\nu,\theta;r)&:=\int_{\R_+} J(x,\nu,\theta;r) \mu_0(dx) \quad \nu \in \mathcal{A}_d, \ r>0.
\end{align}
\begin{remark}[On the initial distribution] We point out that all the results in the previous sections hold true also for profit functionals as in \eqref{def profit functional MFG random initial condition}; that is, if the deterministic initial condition $X_{0^-}=x \in \R_+$ is replaced by $X_{0^-}=\xi$, for a positive square-integrable $\mathcal{F}_0$-random variable $\xi$, with distribution $\mu_0$. 
	In particular, by the Markov property of the solution to the reflected Skorokhod problem (cf.\ Theorem 1.2.2 and Exercise 1.2.2 in \cite{Pilipenko2014}), the MFG equilibria $(\nu^{\star,r},\theta_r^\star)$ and $(\nu^{\star,e},\theta_e^\star)$ are still characterized by couples $(x_r^\star,\theta_r^\star)$ and $(x_e^\star,\theta_e^\star)$ solving the systems of equations provided in Theorems \ref{thm discMFG-ex} and \ref{thm ergMFG-ex}, respectively. 
\end{remark}

For any $i=1,\dots, N$, let us consider the policy $\bar \nu^{i,e} \in {\mathcal{A}}_e$ according to which the state is reflected upward at the boundary $x_e^\star$. Similarly, for $r>0$, the policy $\bar \nu^{i,r}\in {\mathcal{A}}_d$ makes the state upward reflected at $x_r^{\star}$. 
We observe that, for $i=1,...,N$ and $q\in\{d,e\}$, the profile strategies $(\bar \nu^{1,q},\dots,\bar \nu^{i-1, q},\bar \nu^{i+1,q},\dots, \bar \nu^{N,q}) \in \widehat{\mathcal{A}}^{N-1}_q$. 
Then, define accordingly:
\begin{equation}
	\label{eq strategies}
	\bar{\boldsymbol{\nu}}^e:=(\bar\nu^{1,e},\dots,\bar\nu^{N,e}), 
	\quad 
	\bar{\boldsymbol{\nu}}^r:=(\bar\nu^{1,r},\dots,\bar\nu^{N,r}), \quad  \theta_e^{i,N}:=\theta^N_{\bar{\boldsymbol{\nu}}^{-i,e}}, \quad \text{and} \quad \theta_r^{i,N}:=\theta^N_{\bar{\boldsymbol{\nu}}^{-i,r}}.
\end{equation}

To facilitate our discussion, we enforce some additional mild requirements on the dynamics of the state processes and on the profit function.
\newpage

\begin{assumption}
	\label{ass pi-theta additional}  \ 
	\begin{enumerate}
		\item There exists $x_{b,\sigma}>0$ such that $2\, x\, b(x) + \sigma ^2 (x) \leq 0$ for any $x \geq x_{b,\sigma}$; 
		\item For any $a>0$, there exists a constant $C>0$ such that 
		\[|\pi(x,\theta_1)-\pi(x,\theta_2)|\leq C(1+|x|)|\theta_1-\theta_2|,\quad \forall\theta_1,\theta_2\geq a,\]
		for all $x\in\mathbb R$.
	\end{enumerate}
\end{assumption}
Notice that the previous conditions are satisfied by the benchmark cases in which $b(x)=-\delta x$ or $b(x)=\delta(\lambda-x)$ and $\sigma(x)=\sigma x$ (that is, geometric or affine dynamics) when $2\delta \geq \sigma^2$, and for a profit function $\pi(x,\theta)=x^{\beta}\,\theta^{-(1 + \beta)}$, for some elasticity $\beta\in (0,1)$.

For $\theta>0$ and $r \geq 0$, let  $\widehat{x}_r(\theta)$ be as in Assumption \ref{assumption discounted} and \ref{assumption ergodic}.  
It is easy to show that the function $\widehat{x}_{r}(\theta)$ is continuous in $(\theta,r)$ so that, by the  convergence in Theorem \ref{thm: abelian}, we can set 
$$
\widehat{B}:= 2 \max \left\{ \sup_{r\in (0,1]} \widehat{x}_r(F(f(x_r^\star))), \, \widehat{x}_0(F(f(x_e^\star))), \sup_{r\in (0,1]} {x}_r^\star, \, x_e^\star \right\} < \infty.
$$

Next, for any $i=1,...,N$, by definition of $\bar{\nu}^{i,e}$, we have $X_t^{i,\nu^{i,e}} \geq x_e^\star, \, \P$-a.s., for any $t>0$. 
This fact, for $\theta_{e}^{i,N}$ as in \eqref{eq strategies}, by monotonicity of $f$ and $F$ implies that $\theta_{e}^{i,N} \geq F(f(x_e^\star))$, $\P$-a.s. 
In the same way, $\theta_{r}^{i,N} \geq F(f(x_r^\star))$, $\P$-a.s.\ for each $i=1,...,N$ and $r>0$.
Therefore, since for $r \geq0$ the functions $\widehat{x}_r$ are nonincreasing in $\theta$, by definition of $\widehat{B}$
we have
\begin{equation}\label{eq estimate barrier}
	\widehat{x}_0(\theta_{e}^{i,N}) \leq \widehat{x}_0(F(f(x_e^\star)))\leq \widehat{B}, 
	\quad   
	\widehat{x}_r(\theta_{r}^{i,N}) \leq \widehat{x}_r(F(f(x_r^\star))) \leq \widehat{B}, \  r>0.
\end{equation}

It is shown in the proofs of Theorems \ref{thm discMFG-ex} and \ref{thm ergMFG-ex} that, for any $\theta>0$, the optimal control never acts when the optimally controlled state process lies in the set $\{ y : \pi_x(y,\theta)-(r-b'(y))<0 \} = \{y: y> \widehat{x}_r(\theta) \}$. This observation, together with \eqref{eq estimate barrier}, suggests to concentrate the attention only to those strategies $\nu \in {\mathcal{A}}_q$, $q \in \{d,e\}$, that do not increase when $X^{i,\nu} \geq \widehat{B}.$ In fact, defining, for $q\in \{d,e\}$, the set $\mathcal{A}_q(\widehat{B}):= \{ \nu \in {\mathcal{A}}_q : \, \supp( d\nu) \cap \{X^{i,\nu} \geq \widehat{B} \} = \emptyset  , \P\text{-a.s.} \}$, it can be shown that
\begin{equation}
	\label{eq sup = sup}
	\sup_{\nu \in \mathcal{A}_e} G^i(\nu,\bar{\boldsymbol{\nu}}^{-i,e}) 
	= \sup_{\nu \in \mathcal{A}_e(\widehat{B})} G^i(\nu,\bar{\boldsymbol{\nu}}^{-i,e}), \quad 
	\sup_{\nu \in \mathcal{A}_d} J^i(\nu,\bar{\boldsymbol{\nu}}^{-i,r};r) = \sup_{\nu \in \mathcal{A}_d(\widehat{B})} J^i(\nu,\bar{\boldsymbol{\nu}}^{-i,r};r),
	\ \text{for } r>0.
\end{equation}
Moreover, since $x_r^\star, x_e^\star \leq \widehat{B}/2$, we have  $\bar{\nu}^{i,r}\in \mathcal{A}_d(\widehat{B})$ and $\bar{\nu}^{i,e}\in \mathcal{A}_e(\widehat{B})$. 

All these facts allow to prove the following important a priori estimates. Their proof is in Appendix \ref{proof:lemma a priori estimates}.

\begin{lemma}
	\label{lemma a priori estimates} We have  
	$$
	\sup_{\nu \in \mathcal{A}_e(\widehat{B})} \limsup_{T\to \infty } \frac{1}{T} \E\bigg[ \int_0^T |X_t^{i,\nu}|^2 dt \bigg] < \infty 
	\quad \text{and} \quad
	\sup_{\nu \in \mathcal{A}_d(\widehat{B})}  \E\bigg[ \int_0^T e^{-rt} |X_t^{i,\nu}|^2 dt \bigg] < \infty \ \text{for } r>0.
	$$
\end{lemma}  

We are finally ready to state (and prove in Appendix \ref{proof:thm approximation N-player game}) the main result of this section. It states that mean field solutions realize approximate Nash equilibria in the related symmetric $N$-player games defined in Definition \ref{def:epsilonNash}, when $N$ is large and/or $r$ is small.

\begin{theorem}
	\label{thm approximation N-player game}
	The following approximations hold true:
	\begin{enumerate}
		\item $\bar{\boldsymbol{\nu}}^e$ is an $\varepsilon_N$-NE for the ergodic $N$-player game with $\varepsilon_N \to 0$ as $N\to \infty$;
		\item $\bar{\boldsymbol{\nu}}^r$  is an $\varepsilon_{N,r}$-NE for the ergodic $N$-player game with  $\varepsilon_{N,r} \to 0$ as $N\to \infty$ and $r\to 0$;
		\item $\bar{\boldsymbol{\nu}}^r$ is an $\varepsilon_N$-NE for the discounted $N$-player game with $\varepsilon_N \to 0$ as $N\to \infty$;
		\item $\bar{\boldsymbol{\nu}}^e$ is an $\varepsilon_{N,r}$-NE for the discounted $N$-player game with  $\varepsilon_{N,r} \to 0$ as $N\to \infty$ and $r\to 0$.
	\end{enumerate}
\end{theorem}

\begin{remark}
	We point out that results analogous to Claims 1 and 2 in Theorem \ref{thm approximation N-player game} can be obtained even if the mean field interaction term \eqref{eq definition of theta N} in the $N$-player game is replaced by a time-dependent interaction. As a matter of fact, one can consider 
	$$
	\Theta^{N}_{\boldsymbol{\nu}^{-i}}(t):= F\Big( \frac{1}{N-1} \sum_{j \ne i} f(X_t^{j, \nu^j}) \Big),
	$$
	and define accordingly, for $\nu^i \in \mathcal{A}_e$, player $i$'s ergodic profit functional
	\begin{equation*} 
		G^i(\nu^i,\boldsymbol{\nu}^{-i}):=\limsup_{T\to\infty}\frac{1}{T}\E\biggl[\int_0^T\pi\left(X^{i,\nu^i}_{t},\Theta^{N}_{\boldsymbol{\nu}^{-i}}(t)\right) dt-\nu^i_T \biggl], \quad \boldsymbol{\nu}^{-i} \in {\mathcal{A}}_e^{N-1}.
	\end{equation*} 
\end{remark}


\section{Explicit and numerical illustrations in a case study}
\label{sec:casestudy}

In this section, we illustrate the results of the previous sections in a mean field version of a dynamic game of productivity expansion. The productivity of the representative company evolves stochastically as
\begin{equation}
\label{eq: case-dyn}
dX^{\nu}_t= -\delta X^{\nu}_{t} dt + \sigma X^{\nu}_{t} dW_t + d\nu_t, \quad X_{0-} =x >0.
\end{equation}
In \eqref{eq: case-dyn}, the parameter $\delta>0$ measures the natural decay rate of productivity, e.g.\ because of production machines' deterioration; the Brownian motion $W$ models any exogenous shock, such as new technological achievements, that affects productivity level with volatility $\sigma>0$; $\nu_t$ gives the cumulative amount of investment in production or in R\&D made up to time $t$, and which instantaneously makes the company's productivity increase. 

The representative firm interacts with the continuum of other symmetric companies via the operating profit function, which we take of the form
$$\pi(x,\theta):=x^{\beta}\,\theta^{-(1 + \beta)},\quad \theta:=\bigg(\int_{\mathbb{R}_+}x^{\beta}\, \P_{X^{\nu}_{\infty}}(d x)\bigg)^{\frac{1}{\beta}},  \quad \beta \in (0,1).$$ 
As already noticed in Remark \ref{rem:ISO}, such a choice of $\pi$ can be obtained from an isoelastic demand function obtained from Spence-Dixit-Stiglitz preferences (cf.\ \cite{DixitStigl} and footnote 5 in \cite{Achdouetal}). Here, $\theta^{-1}$ can be interpreted as the price index, whose equilibrium level is then determined as a result of the mean field interaction.

\subsection{The solution to the discounted MFG}
\label{sec:casestudy-disc-MFG}

Within this setting Assumptions \ref{A1}, \ref{ass:X-rec}, \ref{A2}, and \ref{assumption discounted} are easily verified, and, employing Theorem \ref{thm discMFG-ex}, the solution to the mean field stationary discounted game as in Definition \ref{def equilibrium discounted} can be explicitly calculated. Namely, one has
$$\theta^{\star}_r =  \big(\rho^{\star}_r\big)^{-\frac{1}{1+\beta}},$$
where
\begin{equation}
    \label{eq: case-disc-mean}
    \rho^{\star}_r:=\left(\frac{1+\frac{2\delta}{\sigma^2}-\beta}{1+\frac{2\delta}{\sigma^2}}\right)^{\frac{1-\beta^2}{2\beta}}\left\{\frac{\sigma^2[n(r)-\beta][1-m(r)]}{2\beta}\right\}^{\frac{1+\beta}{2}},
\end{equation}
with 
\begin{equation*}
    \label{eq: case-disc-m-H}
    \begin{aligned}
    &m(r):=\left(\frac{\delta}{\sigma^2}+\frac{1}{2}\right)-\sqrt{\left(\frac{\delta}{\sigma^2}+\frac{1}{2}\right)^2+\frac{2r}{\sigma^2}}<0,\,\,n(r):=\left(\frac{\delta}{\sigma^2}+\frac{1}{2}\right)+\sqrt{\left(\frac{\delta}{\sigma^2}+\frac{1}{2}\right)^2+\frac{2r}{\sigma^2}}>1.\\
    \end{aligned}
\end{equation*}
Hence, $(\rho^{\star}_r)^{\frac{1}{1+\beta}}=(\theta^{\star}_r)^{-1}$, and -- up to taking the root of order $1+\beta$ -- $\rho^{\star}_r$ identifies with the equilibrium price index. Also, defining 
$$H(\rho^{\star}_r;r):=\frac{2\rho^{\star}_r}{\sigma^2[n(r)-\beta][\beta-m(r)]},$$
the reflection boundary is
\begin{equation}
    \label{eq: case-disc-l}
    x^{\star}_r=\left[\frac{H(\rho^{\star}_r;r)\beta[\beta-m(r)]}{1-m(r)}\right]^{\frac{1}{1-\beta}}.
\end{equation}
Furthermore, the equilibrium stationary distribution is given by
\begin{equation}
\label{eq:statdistr-disc}
\P_{X^{\nu^{\star,r}}_{\infty}}(d x) = \frac{1+\frac{2\delta}{\sigma^2}}{\big(x^{\star}_r\big)^{-\frac{2\delta}{\sigma^2}-1}} x^{-\frac{2\delta}{\sigma^2}-2} \mathds{1}_{[x^{\star}_r, \infty)}(x) dx.
\end{equation}
Figure \ref{subfig: disc-mf-density} describes the density function when $r=0.5$.


\subsection{The solution to the ergodic MFG}
\label{sec:casestudy-erg-MFG}

\begin{figure}
\centering
\begin{subfigure}[t]{0.45\textwidth}
\includegraphics[width = \textwidth]{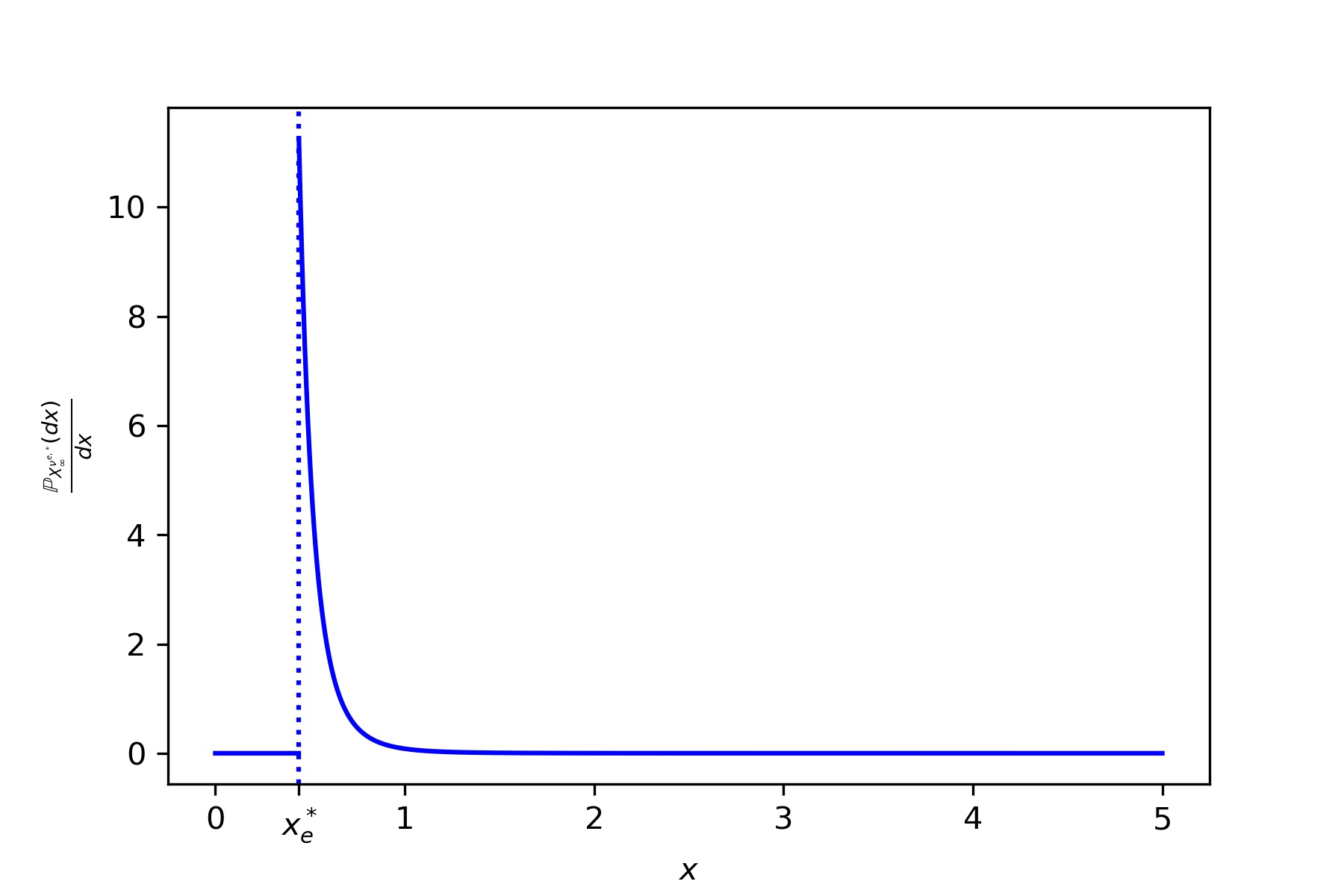}
\caption{Ergodic MFG}
\label{subfig: erg-mf-density}
\end{subfigure}
\begin{subfigure}[t]{0.45\textwidth}
\includegraphics[width = \textwidth]{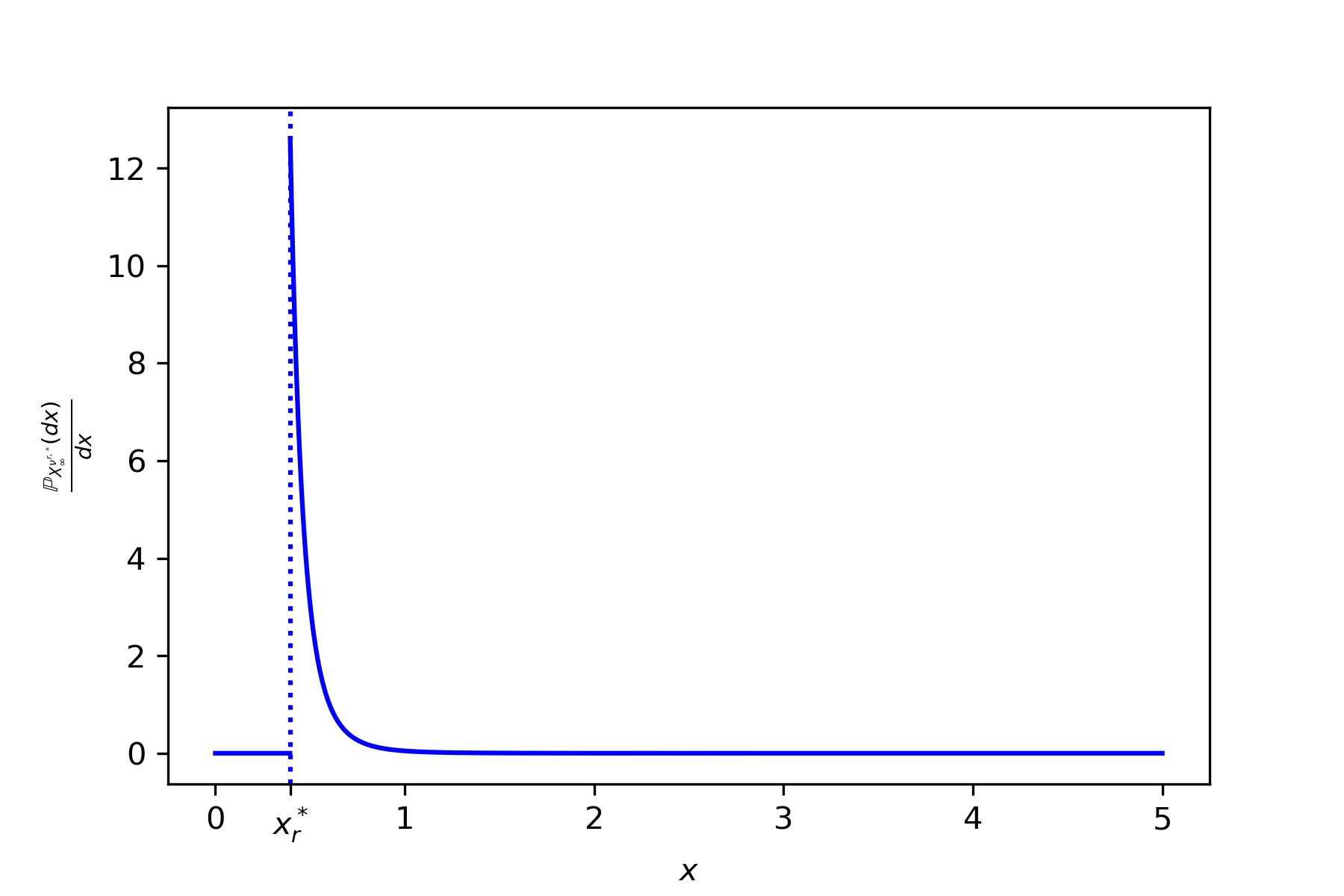}
\caption{Discounted MFG with $r=0.5$}
\label{subfig: disc-mf-density}
\end{subfigure}
\caption{Density functions of the stationary distributions in equilibrium}
\label{fig: mf-density}
\end{figure}

In the setting of this section it is also clear that Assumptions \ref{A1}, \ref{ass:X-rec}, \ref{A2}, and \ref{assumption ergodic} are verified. Hence, Theorem \ref{thm ergMFG-ex} can be applied in order to determine the solution to the ergodic MFG defined in \ref{def equilibrium discounted}. Also in this case the equilibrium barrier and weighted average can be explicitly determined. In particular, 
$$\theta^{\star}_e =  \big(\rho^{\star}_e\big)^{-\frac{1}{1+\beta}},$$
where
\begin{equation}
    \label{eq: case-erg-mean}
    \rho^{\star}_e:=\left[\frac{n-\beta}{n}\right]^{\frac{1-\beta^2}{2\beta}}\left[\frac{\gamma^2(n-\beta)}{2\beta}\right]^{\frac{1+\beta}{2}},
\end{equation}
with $n:=1+\frac{2\delta}{\sigma^2}$.
Moreover, the reflection barrier is
\begin{equation}
    \label{eq: case-erg-l}
    x^{\star}_e=\left[\frac{2\rho^{\star}_e\beta}{\sigma^2(n-\beta)}\right]^{\frac{1}{1-\beta}},
\end{equation}
while the the equilibrium stationary distribution, as in Figure \ref{subfig: erg-mf-density}, is given by
\begin{equation}
\label{eq:statdistr-erg}
\P_{X^{\nu^{\star,e}}_{\infty}}(d x) = \frac{1+\frac{2\delta}{\sigma^2}}{\big(x^{\star}_e\big)^{-\frac{2\delta}{\sigma^2}-1}}x^{-\frac{2\delta}{\sigma^2}-2} \mathds{1}_{[x^{\star}_e, \infty)}(x) dx.
\end{equation}

\begin{figure}
\begin{center}
\includegraphics[scale=0.65]{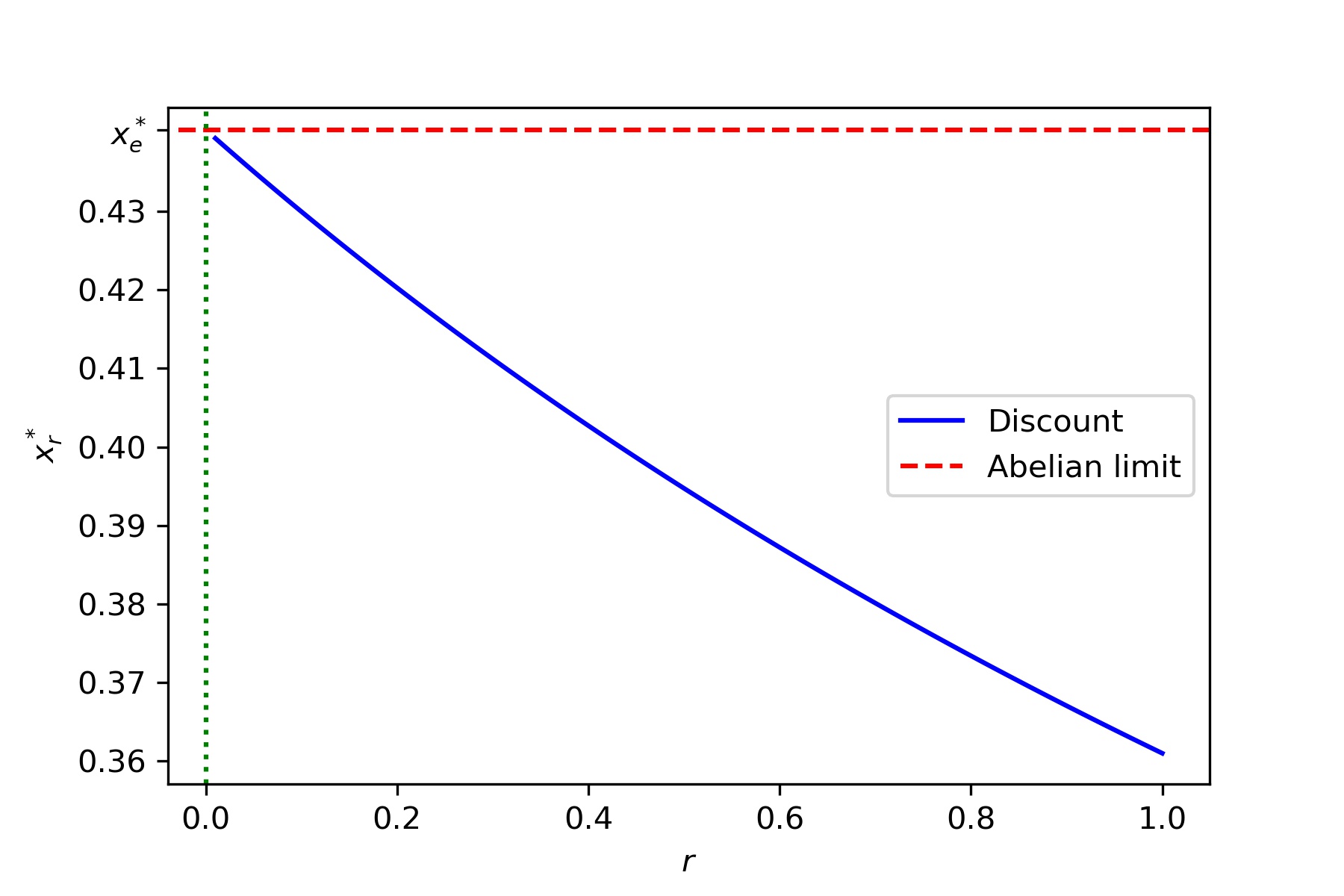}
\caption{Abelian limit: Convergence of $x^{\star}_r$ towards $x^{\star}_e$.}
\label{subfig: abelian-l}
\end{center}
\end{figure}

\begin{figure}
\begin{center}
\includegraphics[scale=0.65]{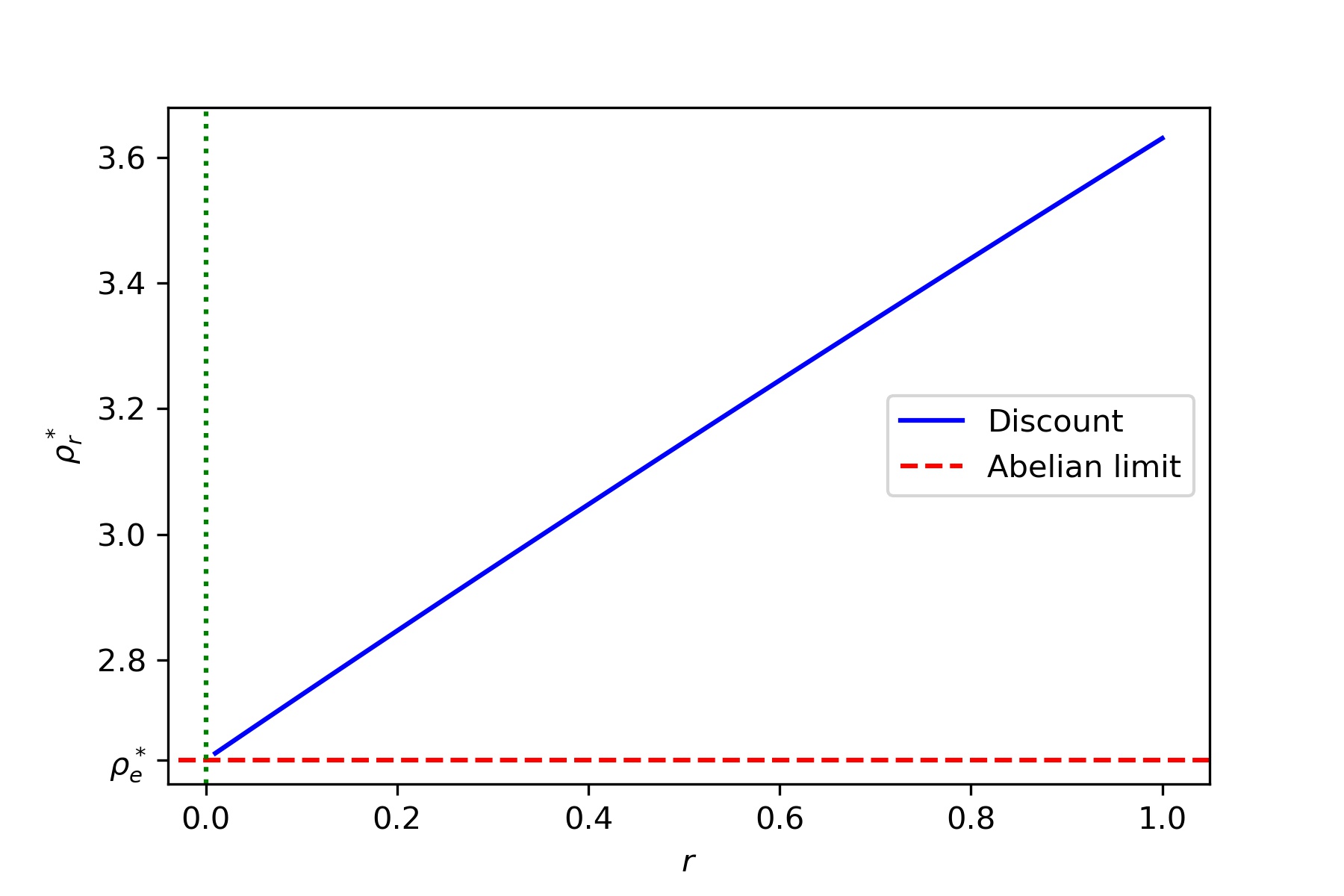}
\caption{Abelian limit: Convergence of $\rho^{\star}_r$ towards $\rho^{\star}_e$.}
\label{subfig: abelian-rho}
\end{center}
\end{figure}

Taking limits in \eqref{eq: case-disc-mean} and \eqref{eq: case-disc-l} as $r\downarrow 0$, by continuity, it is immediate to see that $\theta^{\star}_r \to \theta^{\star}_e$ and $x^{\star}_r \to x^{\star}_e$, consistently to our general result presented in Theorem \ref{thm: abelian}. This is also illustrated in Figures \ref{subfig: abelian-l} and \ref{subfig: abelian-rho} where it is plotted the convergence, with respect to $r$, of $x^{\star}_r$ to $x^{\star}_e$, and of $\rho^{\star}_r$ to $\rho^{\star}_e$, respectively. We see that $x^{\star}_r$ increases as $r$ decreases; i.e.\ a lower discount rate makes the representative company invest earlier. On the other hand, the larger $r$ is, the larger is $\rho^{\star}_r$ since the equilibrium population's productivity is distributed over a bigger interval.

\subsection{Sensitivity analysis of mean field equilibria}
\label{sec:sensitivity}
In this section, we study the sensitivity of the mean-field solutions with respect to key model parameters, namely drift coefficient $\delta$, volatility coefficient $\sigma$, and elasticity $\beta$. The default choices of these parameters are set to be $\delta_{{\rm default}} = 2$, $\sigma_{{\rm default}}=1$ and $\beta_{{\rm default}}=0.6$.

\subsubsection{Drift coefficient $\delta$.}

\begin{figure}
\centering
\begin{subfigure}[t]{0.45\textwidth}
\includegraphics[width=\textwidth]{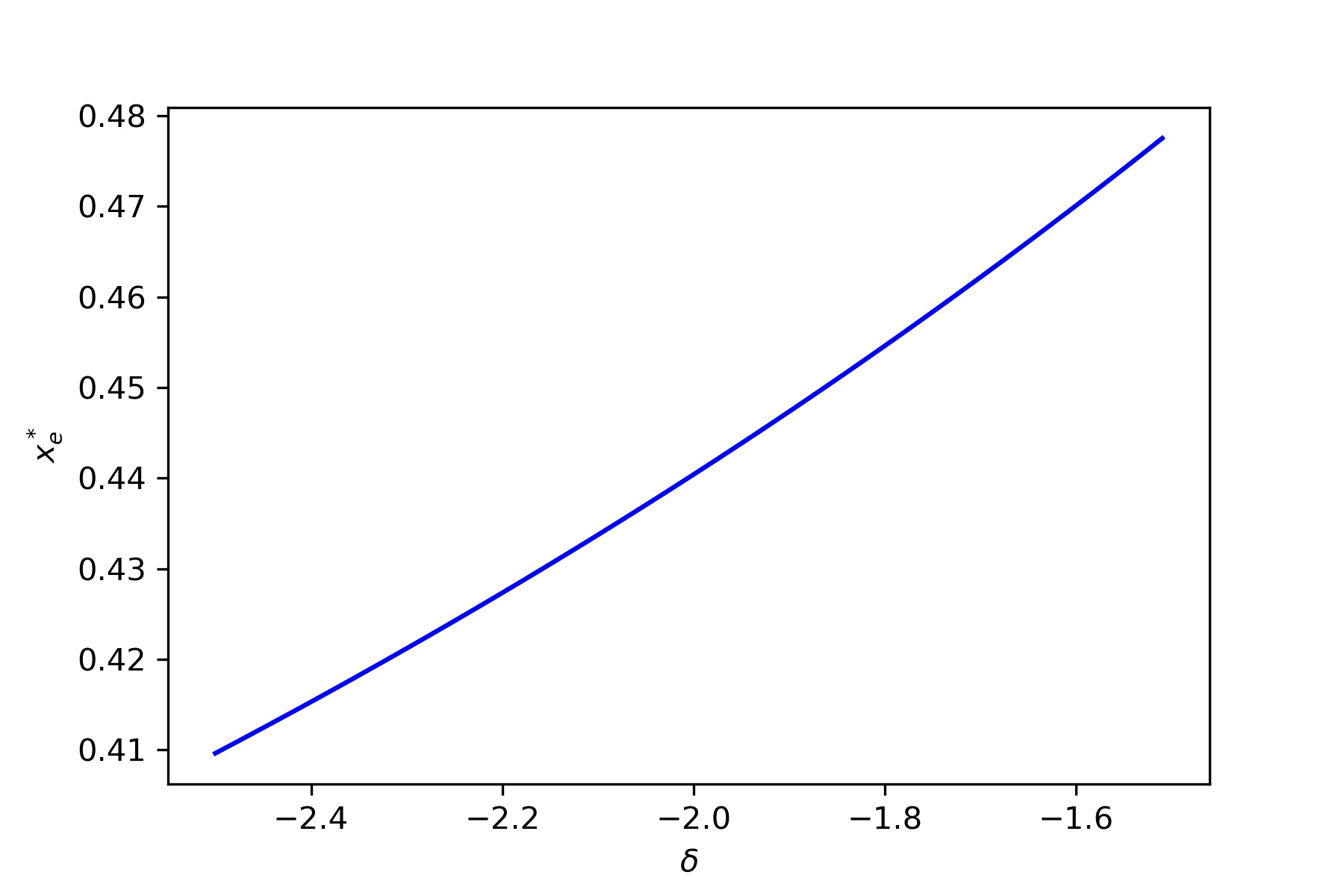}
\caption{Ergodic MFG: $x^{\star}_e$ vs.\ $\delta$}
\label{subfig: erg-l-delta}
\end{subfigure}
\begin{subfigure}[t]{0.45\textwidth}
\includegraphics[width=\textwidth]{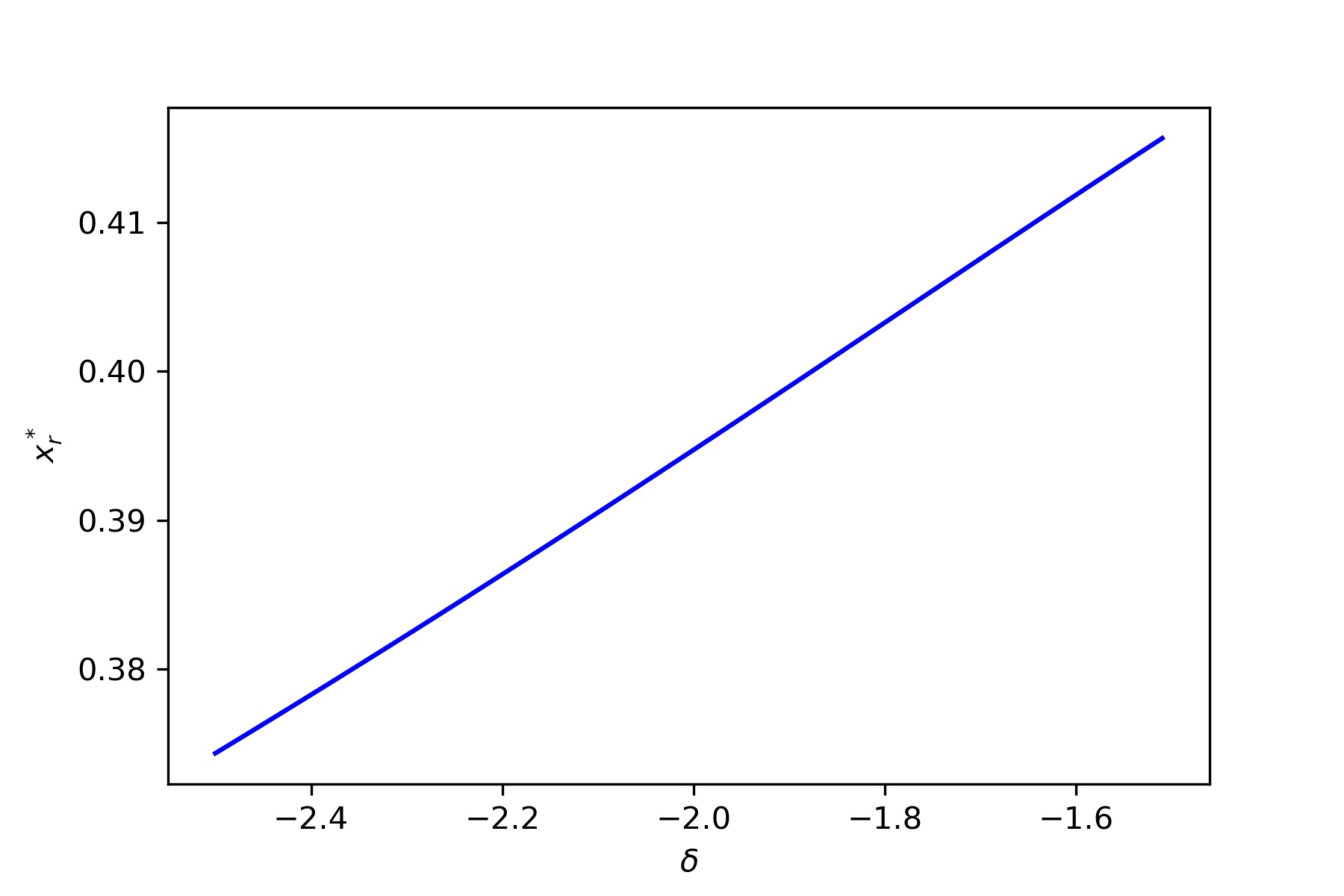}
\caption{Discounted MFG with $r=0.5$: $x^{\star}_r$ vs.\ $\delta$}
\label{subfig: disc-l-delta}
\end{subfigure}
\caption{Sensitivity of reflection boundaries with respect to $\delta$}
\label{fig: l-delta}
\end{figure}

\begin{figure}
\centering
\begin{subfigure}[t]{0.45\textwidth}
\includegraphics[width = \textwidth]{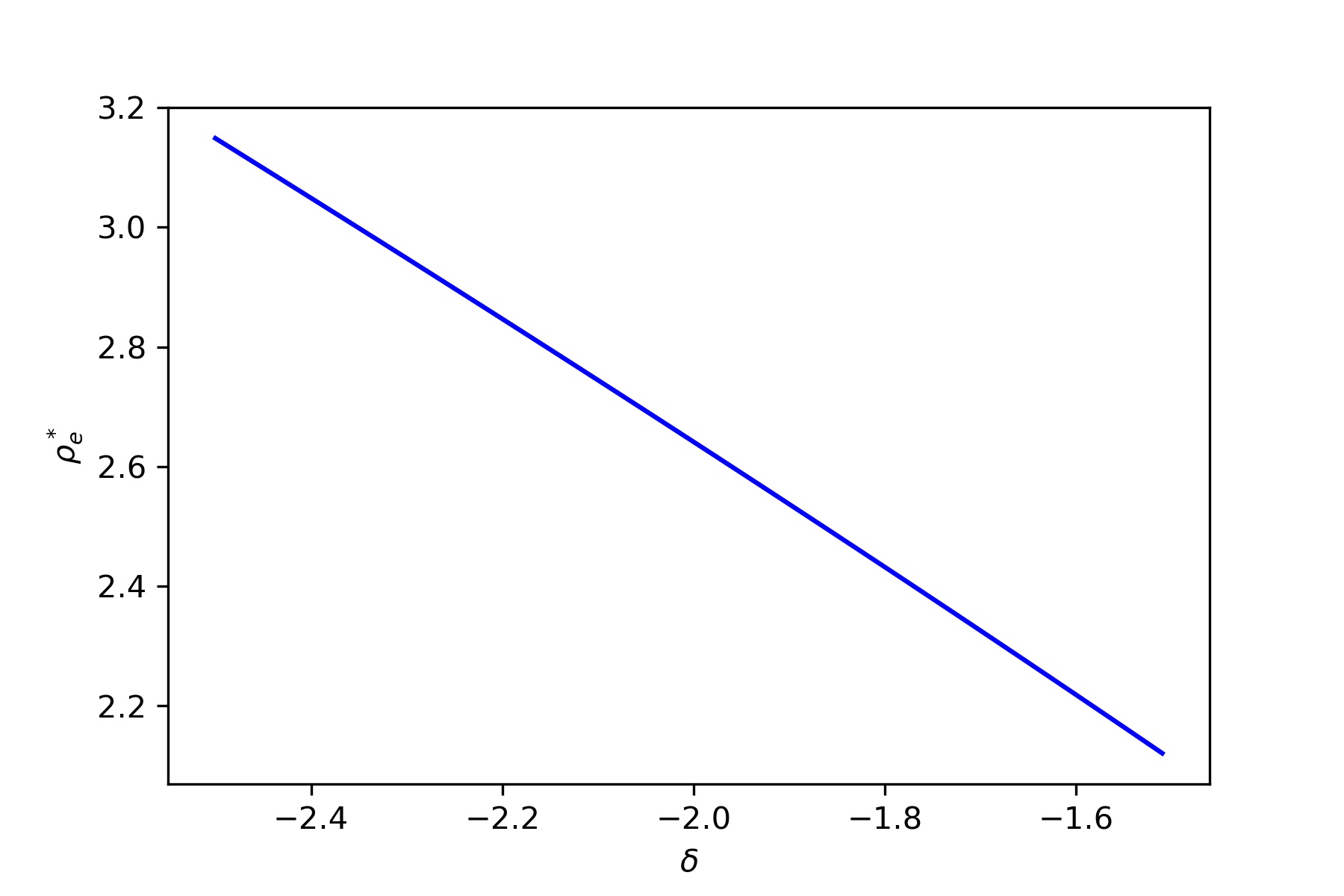}
\caption{Ergodic MFG: $\rho^{\star}_e$ vs.\ $\delta$}
\label{subfig: erg-rho-delta}
\end{subfigure}
\begin{subfigure}[t]{0.45\textwidth}
\includegraphics[width = \textwidth]{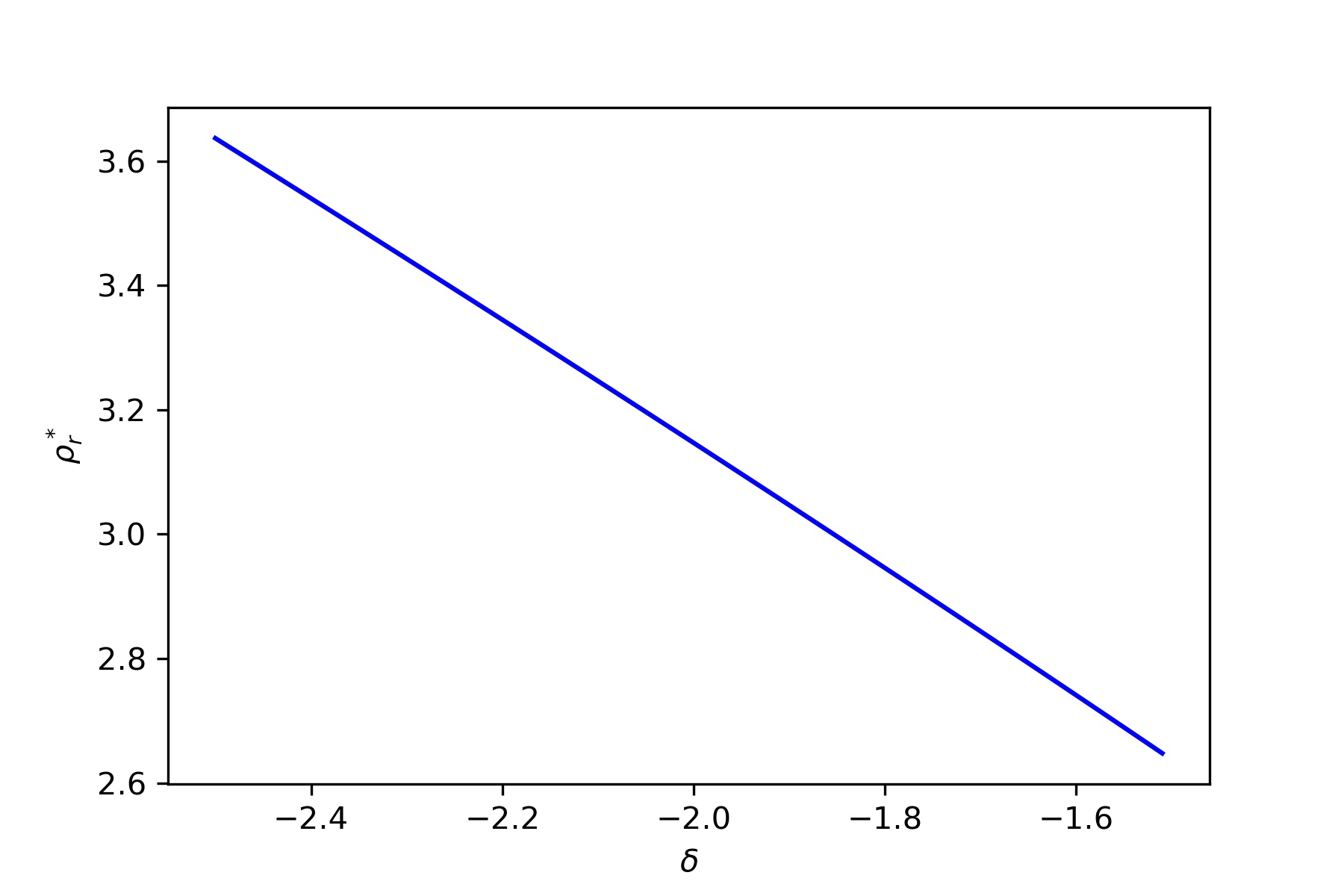}
\caption{Discounted MFG with $r=0.5$: $\rho^{\star}_r$ vs.\ $\delta$}
\label{subfig: disc-rho-delta}
\end{subfigure}
\caption{Sensitivity of $\rho^\star$ with respect to $\delta$}
\label{fig: rho-delta}
\end{figure}

The productivity decreases exponentially (in expectation) at rate $\delta$. The decision then depends on the trade-off between the running payoff $\rho x_t^\beta$, with $\rho:=\theta^{-(1+\beta)}$, and the instantaneous cost of intervention $d\nu_t$. If no intervention occurs during the time interval $[t,t+\Delta t]$, then the running payoff at $t+\Delta t$ is given by
\[\rho x_t^\beta\exp\left\{(-\beta\delta-\frac{\sigma^2}{2}\beta(1-\beta))\Delta t\right\}\exp\left\{\beta\sigma (W_{t+\Delta t}-W_t)-\frac{\beta^2\sigma^2}{2}\Delta t\right\},\]
where $-\beta\delta-\frac{\sigma^2}{2}\beta(1-\beta)$ is the expected growth rate. The larger $\delta$ gets, the larger the rate of depreciation of productivity.
As shown in Figure \ref{subfig: erg-rho-delta}, as larger as $\delta$ gets, indicating a decreasing productivity, the equilibrium price $\rho^{\star}$ increases. Higher price index then leads to higher tolerance of a low individual production level, as shown in Figure \ref{subfig: erg-l-delta}. Similar trends in $\rho_r^*$ and $x_r^*$ as $\delta$ varies can be observed for the discounted mean field game as shown in Figures \ref{subfig: disc-rho-delta} and \ref{subfig: disc-l-delta}.

\subsubsection{Volatility coefficient $\sigma$.}
\begin{figure}
\centering
\begin{subfigure}[t]{0.45\textwidth}
\includegraphics[width = \textwidth]{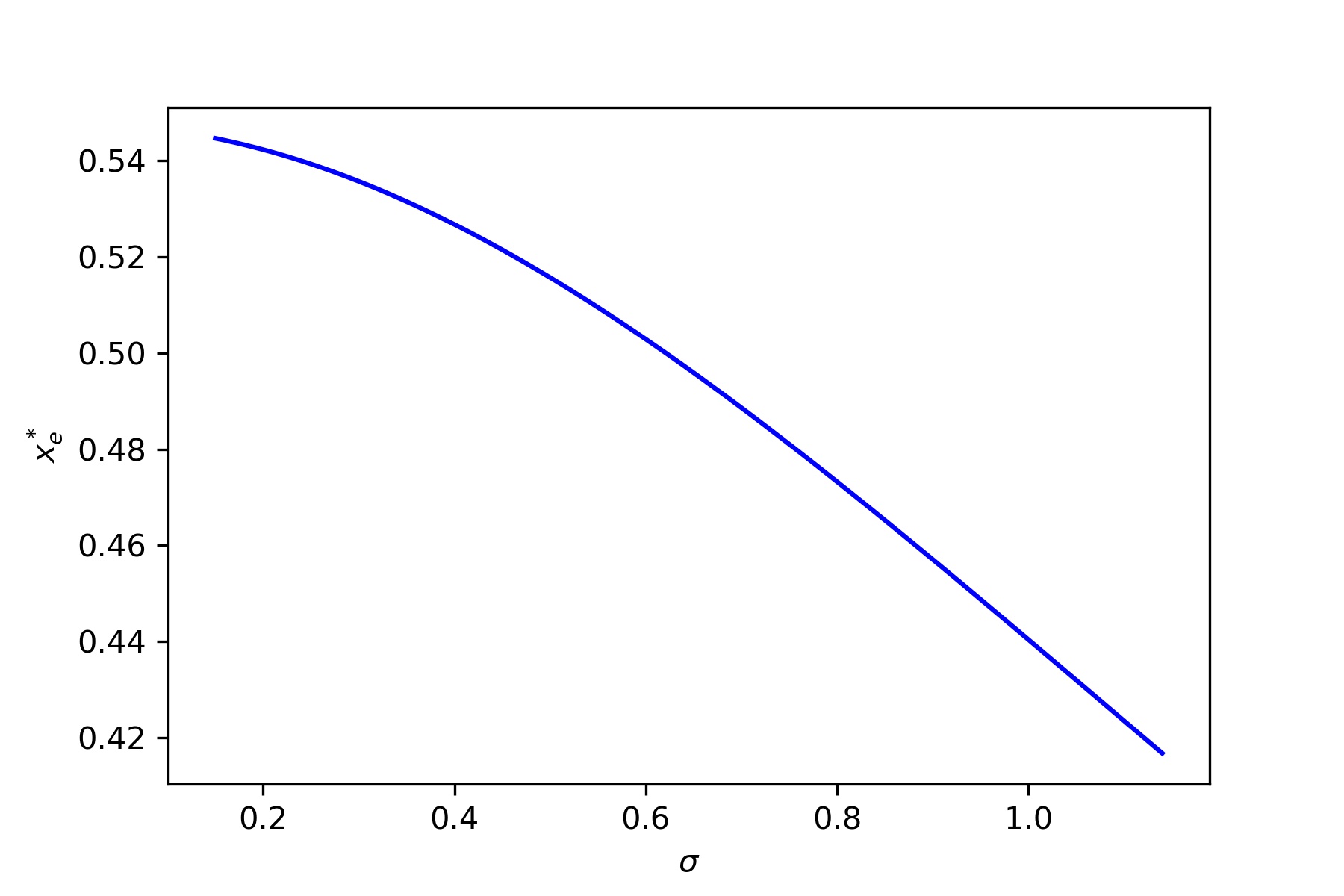}
\caption{Ergodic MFG: $x^{\star}_e$ vs.\ $\sigma$}
\label{subfig: erg-l-gamma}
\end{subfigure}
\begin{subfigure}[t]{0.45\textwidth}
\includegraphics[width = \textwidth]{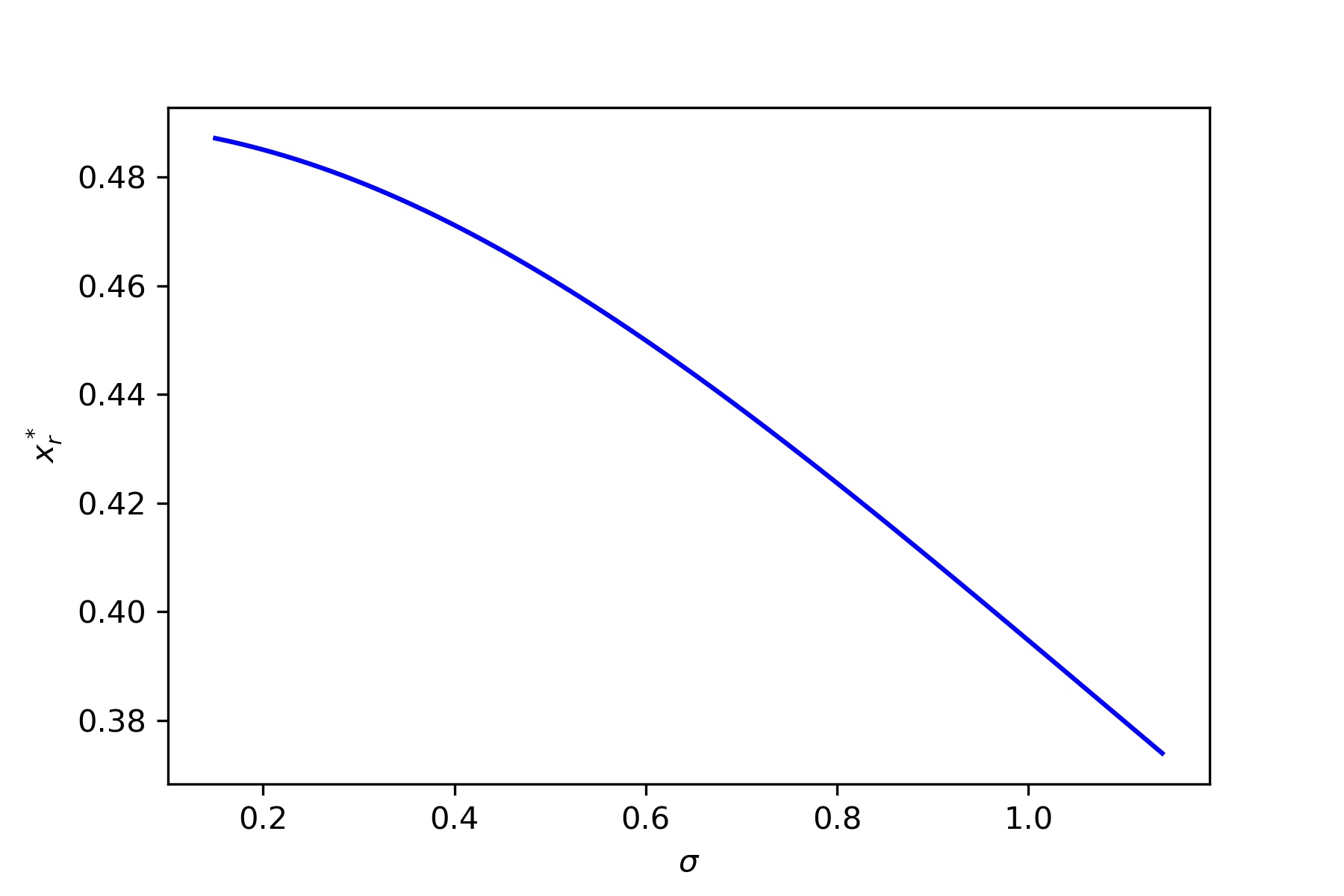}
\caption{Discounted MFG with $r=0.5$: $x^{\star}_r$ vs.\ $\sigma$}
\label{subfig: disc-l-gamma}
\end{subfigure}
\caption{Sensitivity of the reflection boundaries with respect to $\sigma$}
\label{fig: l-gamma}
\end{figure}

\begin{figure}
\centering
\begin{subfigure}[t]{0.45\textwidth}
\includegraphics[width = \textwidth]{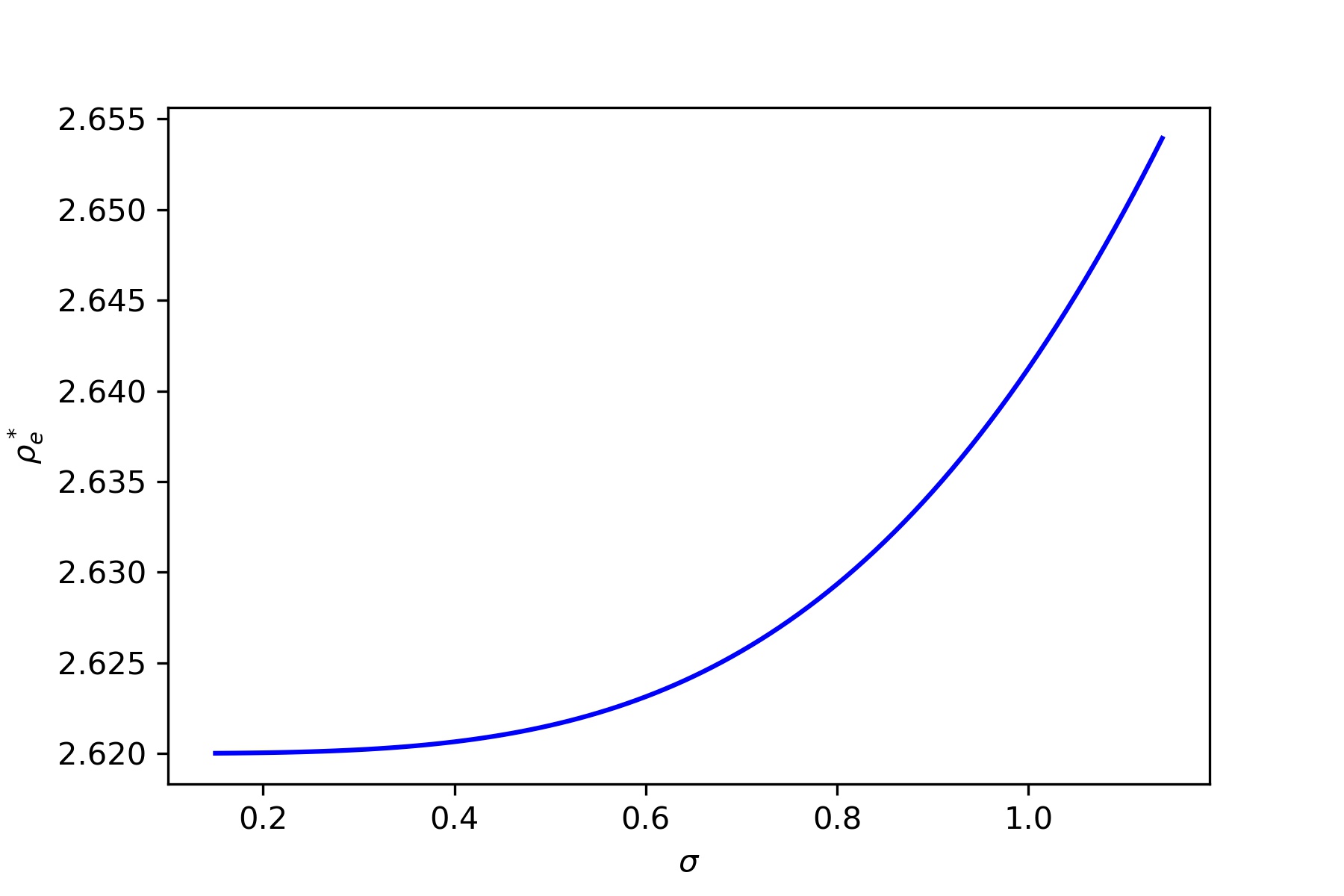}
\caption{Ergodic MFG: $\rho^{\star}_e$ vs.\ $\sigma$}
\label{subfig: erg-rho-gamma}
\end{subfigure}
\begin{subfigure}[t]{0.45\textwidth}
\includegraphics[width = \textwidth]{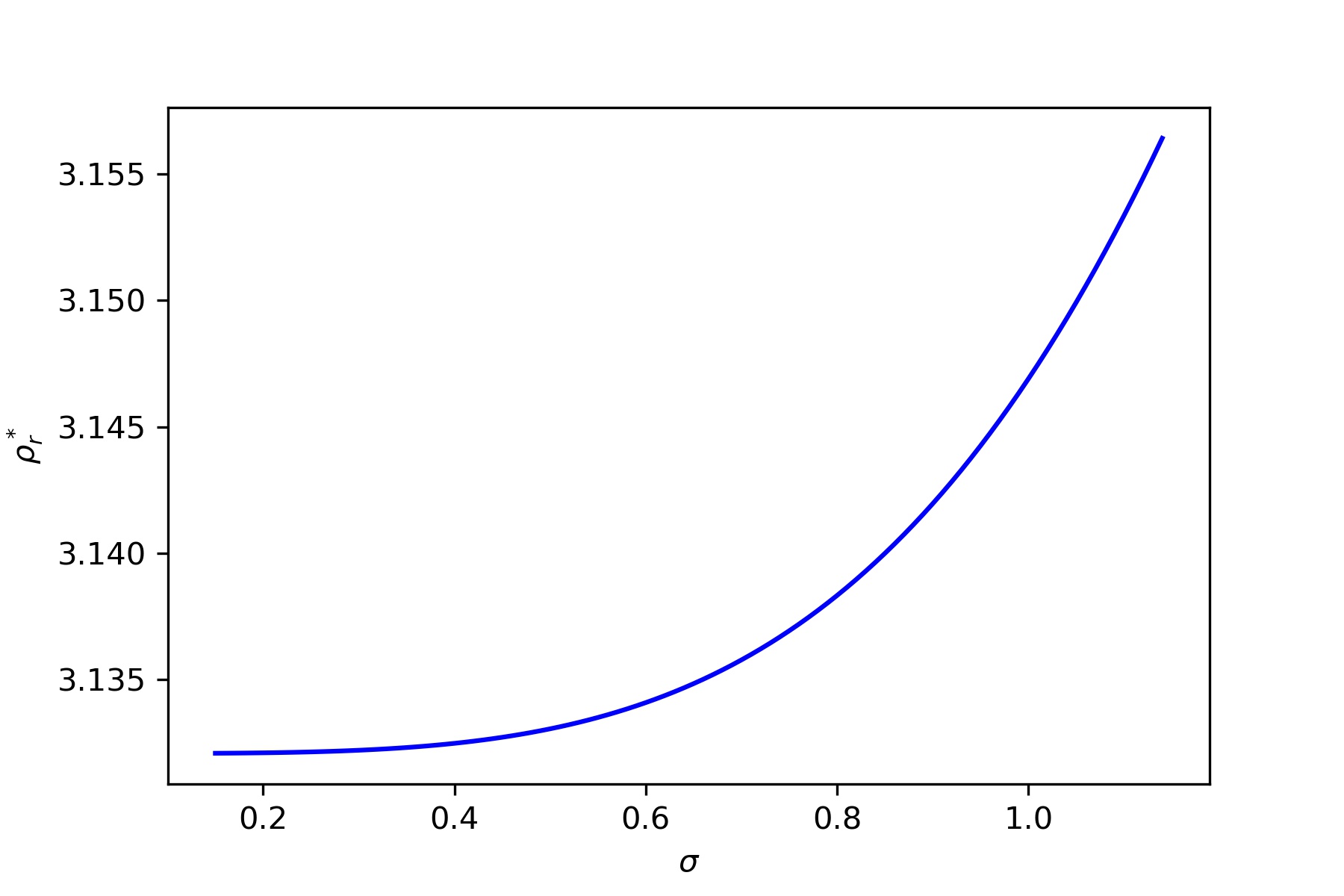}
\caption{Discounted MFG with $r=0.5$: $\rho^{\star}_r$ vs.\ $\sigma$}
\label{subfig: disc-rho-gamma}
\end{subfigure}
\caption{Sensitivity of $\rho^\star$ with respect to $\sigma$}
\label{fig: rho-gamma}
\end{figure}

The volatility coefficient $\sigma$ measures the fluctuations in the productivity level. Due to the relation between equilibrium price $\rho^{\star}$ and the limiting productivity level $X_\infty$ given by \eqref{eq: case-erg-mean}, we can see that higher volatility leads to a higher equilibrium price $\rho^{\star}$, as shown in Figure \ref{subfig: erg-rho-gamma}. Meanwhile, players may take advantage of a higher volatility by postponing their interventions; therefore, ``the value of waiting'' increases and as a consequence on can observe a decease in the reflection boundary as $\sigma$ increases in Figure \ref{subfig: erg-l-gamma}. Similar trends in $\rho_r^{\star}$ and $x_r^{\star}$ as $\sigma$ varies can be observed for the discounted mean field game as shown in Figures \ref{subfig: disc-rho-gamma} and \ref{subfig: disc-l-gamma}.

\subsubsection{Elasticity $\beta$.}
\begin{figure}
\centering
\begin{subfigure}[t]{0.45\textwidth}
\includegraphics[width = \textwidth]{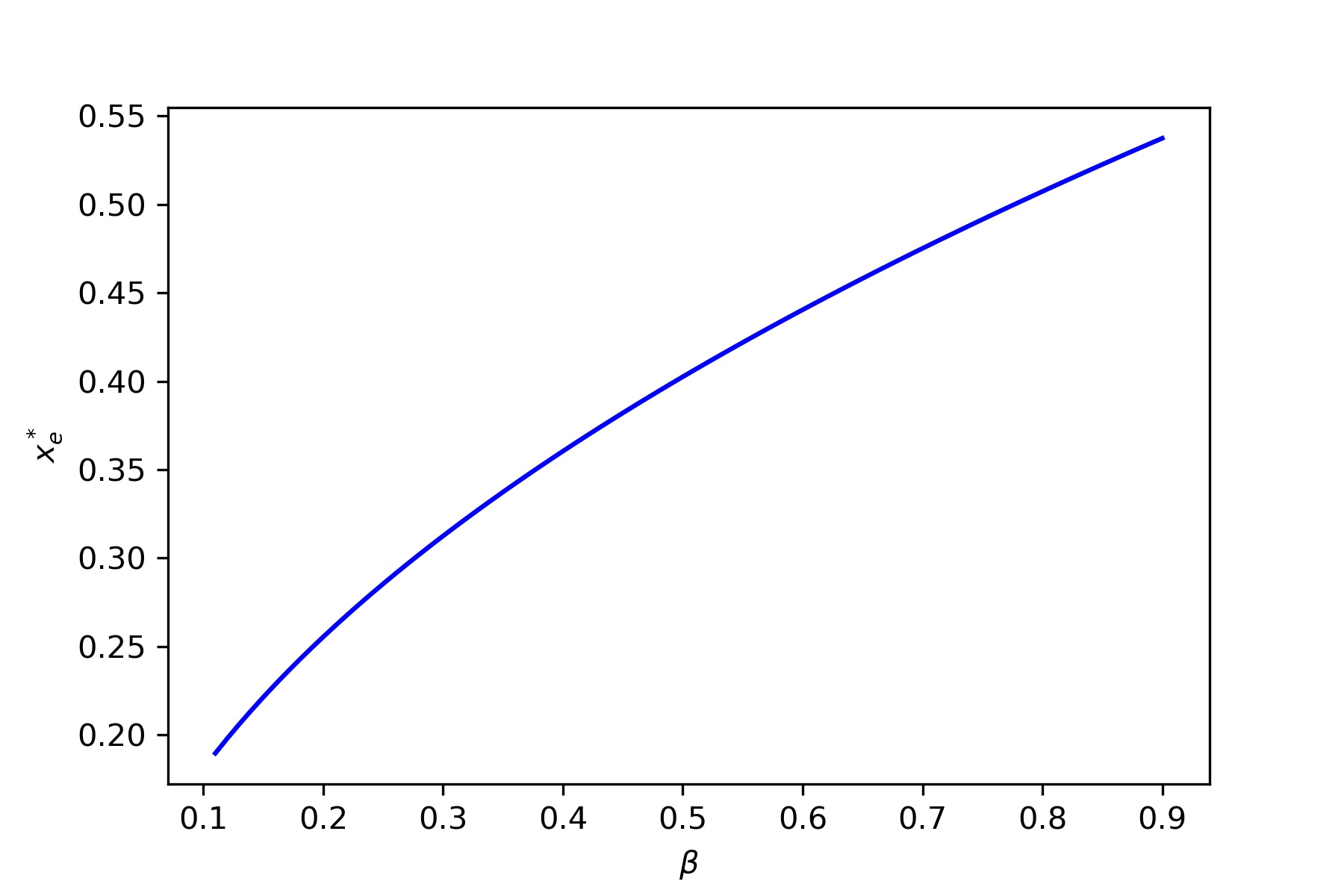}
\caption{Ergodic MFG: $x^{\star}_e$ vs.\ $\beta$}
\label{subfig: erg-l-alpha}
\end{subfigure}
\begin{subfigure}[t]{0.45\textwidth}
\includegraphics[width = \textwidth]{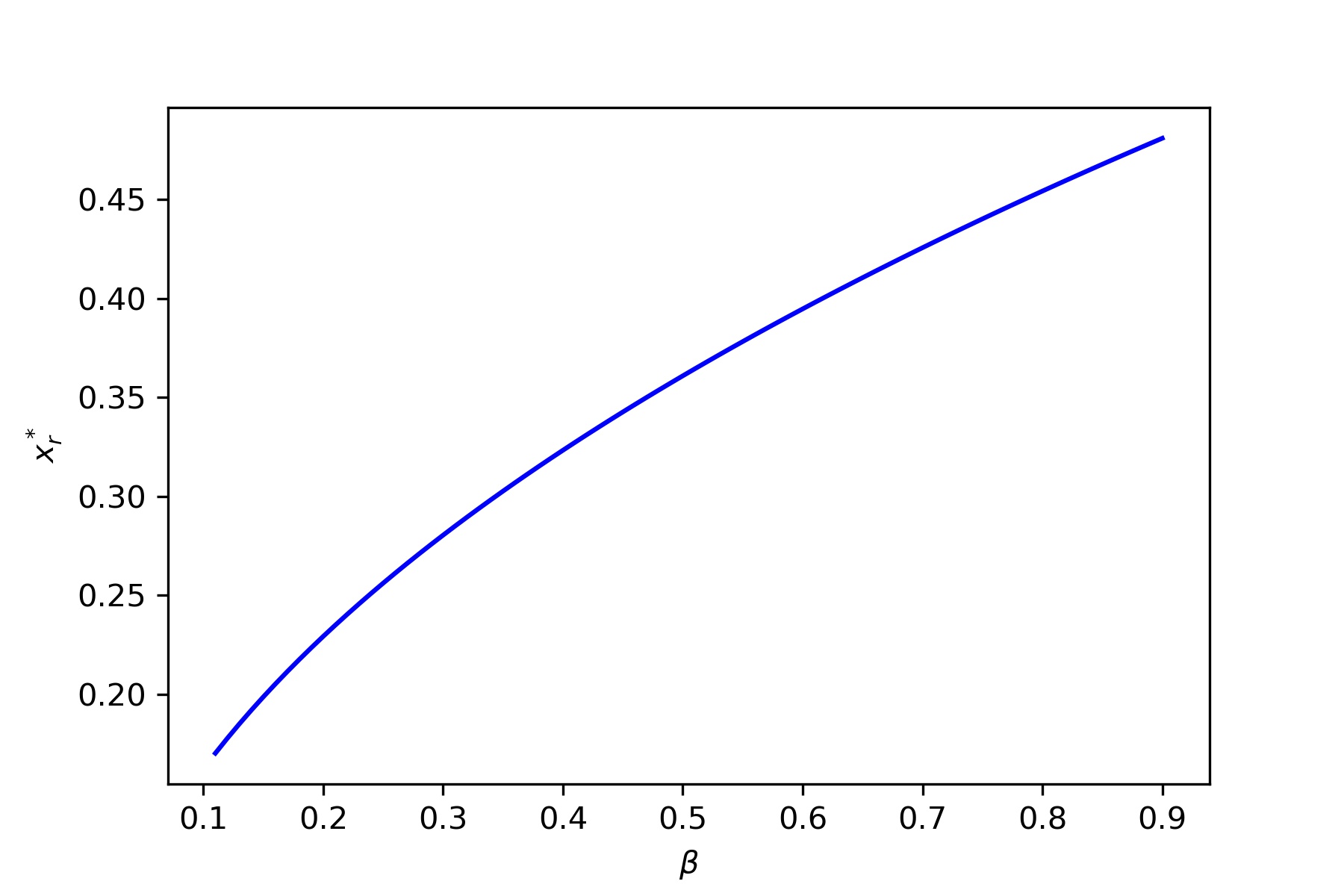}
\caption{Discounted MFG with $r=0.5$: $x^{\star}_r$ vs.\ $\beta$}
\label{subfig: disc-l-alpha}
\end{subfigure}
\caption{Sensitivity of the reflection boundaries with respect to $\beta$}
\label{fig: l-alph}
\end{figure}

\begin{figure}
\centering
\begin{subfigure}[t]{0.45\textwidth}
\includegraphics[width=\textwidth]{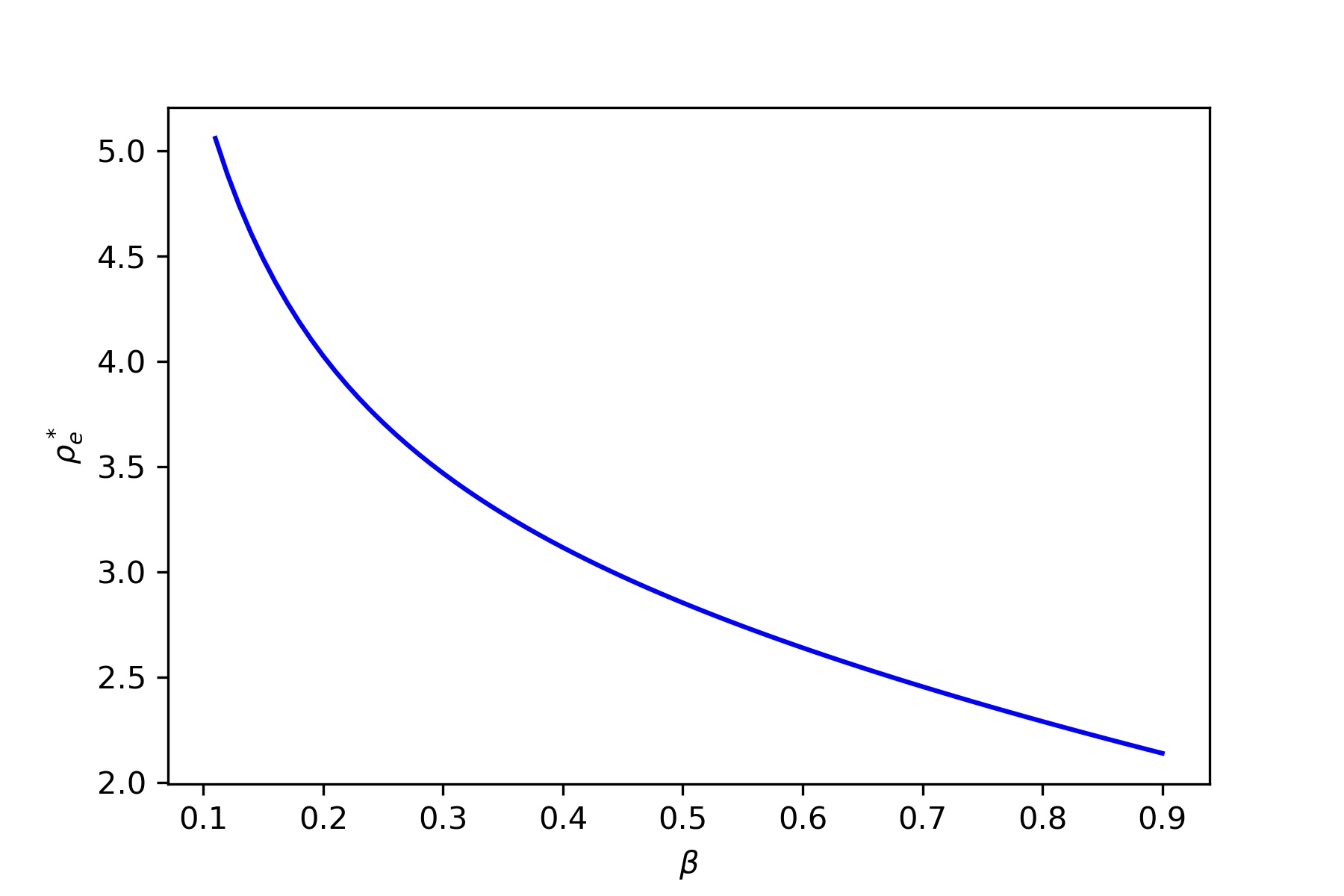}
\caption{Ergodic MFG: $\rho^{\star}_e$ vs.\ $\beta$}
\label{subfig: erg-rho-alpha}
\end{subfigure}
\begin{subfigure}[t]{0.45\textwidth}
\includegraphics[width=\textwidth]{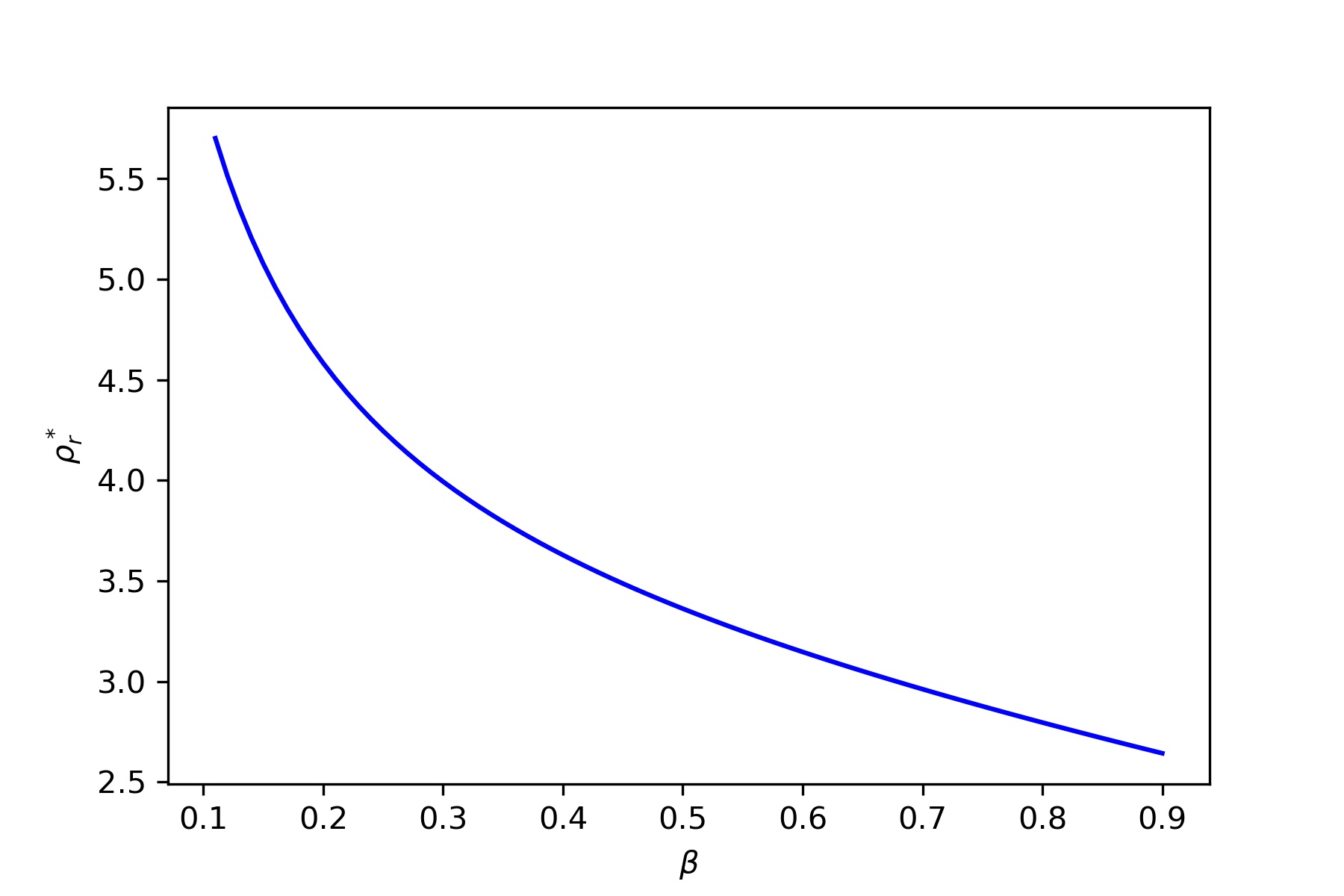}
\caption{Discounted MFG with $r=0.5$: $\rho^{\star}_r$ vs.\ $\beta$}
\label{subfig: disc-rho-alpha}
\end{subfigure}
\caption{Sensitivity of $\rho^\star$ with respect to $\beta$}
\label{fig: rho-alph}
\end{figure}

The parameter $\beta$ provides the elasticity of the payoff with respect to the productivity. The equilibrium price $\rho^{\star}$ inversely depends on $\beta$, showing a decrease as $\beta$ increases in Figure \ref{subfig: erg-rho-alpha}. On the other hand, the more sensitive the profit becomes with respect to product, the more willing the players would like to keep their productivity at a high level. Therefore, as $\beta$ increases, the reflecting boundary also increases, as shown in Figure \ref{subfig: erg-l-alpha}. Similar trends in $\rho_r^{\star}$ and $x_r^{\star}$ as $\beta$ varies can be observed for the discounted mean field game as shown in Figures \ref{subfig: disc-rho-alpha} and \ref{subfig: disc-l-alpha}.

\section{Conclusions}
\label{sec:concl}

In this paper we have considered stationary mean field games in which a representative player can exert a singular control in order to instantaneously increase the level of an underlying It\^o-diffusion process. The considered class of problems well models a mean field version of a symmetric dynamic game of productivity expansion, in which each company interacts with the opponents through its profit, which depends in a decreasing way on a weighted average of the asymptotic productivities. This is possible since the continuum of firms are affected by idiosyncratic shocks and it is assumed that a form of the Law of Large Numbers is valid. In this sense, our games encompass what pointed out by A.K.\ Dixit and R.S.\ Pindyck (cf.\ p.\ 271 in \cite{DixitPindyck}): ``At the industry level, the shocks and responses of firms can aggregate into long-run stationary conditions, so that the industry output and price are nonrandom. However, the equilibrium level of these variables is affected by the parameters of firm-specific uncertainty''. We propose $r$-discounted and ergodic versions of the game. Under appropriate requirements on the data, we show existence and uniqueness of equilibria of barrier type, and, for the first time in the related literature, we connect the equilibria of the two classes of games as the discount rate $r$ vanishes (Abelian limit). Such a link allows us also to shed light on a class of strategic ergodic models that has not yet been investigated, and that can find natural applications in sustainable development and management of public goods. As a matter of fact, we prove that, as $N$ is large and $r$ is small, the mean field equilibrium of the $r$-discounted problem approximates a Nash equilibrium for an $N$-player game in which each exchangeable agent can play singular actions and aims at optimizing an expected ergodic net profit functional. 

There are several directions towards which the present work can be generalized. For example, it would be interesting to investigate analogous results when the representative agent reacts to the current distribution of the population (rather than to the stationary one), and then to study the transition dynamics of the equilibrium. Further, a multi-dimensional underlying state process could be considered in order understand the role of different sources of uncertainty in the mean field equilibrium. These and other generalizations are left for future research.


%

\section{Appendices}
\label{sec:app}
	
	\subsection{Proof of Theorem \ref{thm discMFG-ex}}
	\label{proof:Thm-discMFG}
	
	The proof is organized in two steps.
	\vspace{0.25cm}
	
	\emph{Step 1.} For any fixed $\theta \in \R_+$, here we solve the problem
	\begin{equation}
		\label{SSC-theta}
		V(x,\theta;r):=\sup_{\nu \in \mathcal{A}_d} \E_x\bigg[ \int_0^{\infty} e^{- r s} \pi(X^{\nu}_s, \theta) d  s -  \int_0^{\infty} e^{- r s} d \nu_s \bigg].
	\end{equation} 
	We shall see that an optimal control for \eqref{SSC-theta} is such that to keep (with minimal effort) the state process above a trigger $x^{\star}(\theta)$. Although the arguments of this step are somehow classical (see, e.g., \cite{JackJonhnsonZervos}) we sketch here their main ideas for the sake of completeness. In the following, in order to simplify exposition, we do not explicitly stress the dependency on $r$, unless strictly necessary.
	
	Motivated by the intuition that a costly investment should be made only when the productivity is sufficiently low, for any $x \in \mathbb{R}_+$ we define the candidate value
	\begin{equation}
		\label{candidatev}
		v(x,\theta):=
		\begin{cases}
			A\phi_r(x) + \overline{v}(x,\theta), & x > x^{\star}(\theta), \\
			(x-x^{\star}(\theta)) + v(x^{\star}(\theta),\theta), & x \leq x^{\star}(\theta),
		\end{cases}
	\end{equation}
	for constants $A$ and $x^{\star}(\theta)$ to be found, and with 
	$$\overline{v}(x,\theta) := \E_x\bigg[ \int_0^{\infty} e^{- r s} \pi(X_s, \theta) d  s\bigg],$$
	which is finite due to $(i)$ and $(ii)$ of Assumption \ref{assumption discounted}.
	
	In order to determine $A$ and $x^{\star}(\theta)$ we impose that $v(\cdot,\theta)$ belongs to $C^2(\mathbb{R}_+)$, from which we obtain that
	\begin{equation}
		\label{eq:B}
		A=-\frac{\overline{v}_{xx}(x^{\star}(\theta),\theta)}{\phi_r''(x^{\star}(\theta))}
	\end{equation}
	and
	$$\overline{v}_x(x^{\star}(\theta),\theta)\phi_r''(x^{\star}(\theta)) - \overline{v}_{xx}(x^{\star}(\theta),\theta)\phi_r'(x^{\star}(\theta)) = \phi_r''(x^{\star}(\theta)).$$
	Now, using that $\phi_r'(x) = - \widehat{\phi}_r(x)$ and dividing both members of the latter by $\S'(x^{\star}(\theta))$ we obtain
	\begin{equation}
		\label{eq:eqxstar-1}
		\frac{\overline{v}_{xx}(x^{\star}(\theta),\theta)\widehat{\phi}_r(x^{\star}(\theta)) - \overline{v}_x(x^{\star}(\theta),\theta)\widehat{\phi}_r'(x^{\star}(\theta))}{\S'(x^{\star}(\theta))}= -\frac{\widehat{\phi}_r'(x^{\star}(\theta))}{\S'(x^{\star}(\theta))}.
	\end{equation} 
	
	Notice now that for any function $h \in C^2(\mathbb{R}_+)$, standard differentiation, and the fact that $\mathcal{L}_{\X}\S =0$ and $(\mathcal{L}_{\X}-(r-b'))g = 0$ for $g \in \{\widehat{\psi},\widehat{\phi}\}$, yield
	\begin{equation}
		\label{derivative1}
		\frac{d}{dx} \left[ \frac{ h'(x)}{\S'(x)} \widehat \phi_r(x) - \frac{\widehat \phi_r'(x)}{\S'(x)} h(x) \right] = \widehat \phi_r(x) \widehat m'(x) \big(\mathcal{L}_{\X} - (r - b'(x))\big) h(x).
	\end{equation}
	This last relation applied to the left-hand side of \eqref{eq:eqxstar-1} with $h=\overline{v}_{x}$, and to the right-hand side of \eqref{eq:eqxstar-1} with $h=1$ gives
	\begin{equation}
		\label{eq:eqxstar-2}
		-\int_{x^{\star}(\theta)}^{\infty} \widehat \phi_r(y) \widehat m'(y) (\mathcal{L}_{\X} - (r - b'(y)))\overline{v}_{x}(y,\theta) d  y =  \int_{x^{\star}(\theta)}^{\infty} \widehat \phi_r(y) (r - b'(y)) \widehat m'(y) d  y.
	\end{equation}
	Using now that $(\mathcal{L}_{\X} - (r - b'(y)))\overline{v}_{x}(y,\theta) = - \pi_x(y,\theta)$ we obtain from \eqref{eq:eqxstar-2} an integral equation for $x^{\star}(\theta)$:
	\begin{equation}
		\label{eq:eqxstar}
		K(x^{\star}(\theta),\theta) = 0, \quad \text{where} \quad K(x,\theta):= \int_{x}^{\infty} \widehat \phi_r(y) \big(\pi_x(y,\theta) - r + b'(y)\big) \widehat m'(y) d y.
	\end{equation}
	Due to Assumption \ref{A2} it is easy to see that $K(\widehat{x}_r(\theta),\theta) < 0$. Moreover,  
	\begin{equation}
		\label{deriv-K}
		K_x(x,\theta) = - \widehat \phi_r(x) \big(\pi_x(x,\theta) - r + b'(x)\big) \widehat m'(x) 
		\left\{ \begin{array}{ll}
			\geq 0, \qquad x \geq \widehat{x}_r(\theta) \\
			< 0, \qquad x < \widehat{x}_r(\theta).
		\end{array} \right.
	\end{equation}
	Also, for any $x < \widehat{x}_r(\theta) - \varepsilon:=\widehat{x}_{\varepsilon}(\theta)$, for suitable $\varepsilon >0$, and for $z \in (x, \widehat{x}_{\varepsilon}(\theta))$, by the integral mean-value theorem we find 
	\begin{align*}
		K(x,\theta) & = \int_{x}^{\widehat{x}_{\varepsilon}(\theta)} \widehat \phi_r(y) \big(\pi_x(y,\theta) - r + b'(y)\big) \widehat m'(y)  d  y + K(\widehat{x}_{\varepsilon}(\theta),\theta) \nonumber \\
		& = \frac{\pi_x(z,\theta) - r + b'(z)}{r - b'(z)} \Big(\frac{\widehat{\phi}_r'(\widehat{x}_{\varepsilon}(\theta))}{\S'(\widehat{x}_{\varepsilon}(\theta))} - \frac{\widehat{\phi}_r'(x)}{\S'(x)}\Big) + K(\widehat{x}_{\varepsilon}(\theta),\theta),
	\end{align*}
	where \eqref{psiphiproperties3} have been used in the last step. Using now that $\pi_x(z,\theta) - r + b'(z) >0$ by Assumption \ref{A2}, that $r - b'(z) \geq 2c>0$, and \eqref{psiphiproperties2bis} we see that that $\lim_{x \downarrow 0} K(x,\theta)= \infty$. The previous considerations thus lead to the existence of a unique $x^{\star}(\theta) \in (0, \widehat{x}_r(\theta))$ solving \eqref{eq:eqxstar}. For later use, we stress that
	\begin{equation}\label{eq Kx <0}
		K_x(x^\star(\theta), \theta) <0.
	\end{equation}
	
	It can then be checked that $v(x,\theta)$ as in \eqref{candidatev} is a $C^2$-solution to the HJB equation
	\begin{equation}\label{eq HJB discounted}
		\min\big\{(\mathcal{L}_{X} - r)u(x,\theta) + \pi(x,\theta), 1 - u_x(x,\theta)\big\}=0.
	\end{equation}
	In turn, this allows to show, via a classical verification theorem, that $v(x,\theta)=V(x,\theta)$ and that the control $\nu^{\star}(\theta)$ such that 
	\begin{equation}
		\label{eq:Sk}
		X^{x,\nu^{\star}(\theta)}_t \geq x^{\star}(\theta) \quad \text{and} \qquad \nu^{\star}_t(\theta)=\int_{0}^t\mathds{1}_{\{X^{x,\nu^{\star}(\theta)}_s \leq x^{\star}(\theta) \}} d \nu^{\star}_s(\theta), \quad \forall t \geq 0\,\, \P-\text{a.s.},
	\end{equation}
	belongs to $\mathcal{A}_d$ and is optimal. As a matter of fact, since the free boundary $x^{\star}(\theta)$ is a constant, the latter control rule exists by classical results on the Skorokhod reflection problem (cf.\ Chapter 6 in \cite{Harrison} and Chapter 3.6 in \cite{KS})
	\vspace{0.25cm}
	
	\emph{Step 2.} Since the control $\nu^{\star}(\theta)$ reflects upward the process $X^{x,\nu^{\star}(\theta)}$ at $x^{\star}(\theta)$, the optimally controlled process $X^{x,\nu^{\star}(\theta)}$ is positively recurrent and its stationary distribution is such that (cf.\ Section 12 of Chapter II in \cite{Borodin Salminen})
	$$\P_{X^{x,\nu^{\star}(\theta)}_{\infty}}(dx)= \frac{m'(x)\mathds{1}_{[x^{\star}(\theta),\infty)}(x)}{\int_{x^{\star}(\theta)}^{\infty} m'(y) d  y} dx.$$
Indeed, by Assumption \ref{ass:X-rec} we have 
	\begin{equation}
		\label{X-posrec}
		\int_{x^{\star}(\theta)}^{\infty} m'(y) d  y < \infty.
	\end{equation}
	It thus follows from \eqref{X-posrec} that the consistency equation (i.e., (2) in Definition \ref{def equilibrium discounted}) reads
	\begin{equation*}
		\theta = F \bigg( \int_{x^{\star}(\theta)}^{\infty} f(y) \P_{X^{x,\nu^{\star}(\theta)}_{\infty}}(d y) \bigg)=  F \left(\frac{\int_{x^{\star}(\theta)}^{\infty} f(y) m'(y) d  y}{\int_{x^{\star}(\theta)}^\infty m'(y) d  y } \right);
	\end{equation*} 
	that is,
	\begin{equation}
		\label{eq:CE}
		Q(\theta):=\int_{x^{\star}(\theta)}^{\infty} (f(y) - F^{-1}(\theta) ) m'(y) d  y = 0. 
	\end{equation}  
	
	We now show that \eqref{eq:CE} admits a unique solution $\theta^{\star}$ so that $(\nu^{\star}, \theta^{\star}):=(\nu^{\star}(\theta^{\star}), \theta^{\star})$ is the mean field equilibrium for the discounted stationary MFG.
	
	Recall the definition of $x^\star(\theta)$ and $K$ in \eqref{eq:eqxstar}. Since $K \in C^1(\mathbb{R}_+^2)$ due to $\pi$-$(iv)$ in Assumption \ref{A2}, by the implicit function theorem we find that $\theta \mapsto x^{\star}(\theta)$ is continuously differentiable and has derivative 
	\begin{equation}\label{eq dx : d theta} 
		\frac{d }{d \theta}x^{\star}(\theta) = - \frac{K_{\theta}(x^{\star}(\theta),\theta)}{K_{x}(x^{\star}(\theta),\theta)} < 0,
	\end{equation}
	where the last inequality follows from \eqref{eq Kx <0} and from the fact that $\theta \mapsto \pi_x (x,\theta)$ is strictly decreasing, by Assumption \ref{A2}.
	
	We thus have that $Q$ as in \eqref{eq:CE} is continuously differentiable with derivative
	\begin{equation}
		\label{eq:derG}
		\frac{d }{d \theta}Q(\theta) = - (f(x^{\star}(\theta)) - F^{-1}(\theta)) m'(x^{\star}(\theta)) \frac{d }{d \theta}x^{\star}(\theta) - \frac{1}{F'(F^{-1}(\theta))}\int_{x^{\star}(\theta)}^{\infty} m'(y) d  y.
	\end{equation}
	
	Let now $\widehat{\theta}$ be the unique solution to $f(x^{\star}(\theta)) - F^{-1}(\theta) =0$. Such a value indeed exists. To see this notice that $f \circ x^{\star} $ is strictly decreasing and continuous. Moreover, by using $\pi$-$(iii)$ of Assumption \ref{A2} it can be shown that $x^{\star}(\theta) \rightarrow +\infty$ as $\theta \downarrow 0$ and $x^{\star}(\theta) \rightarrow 0$ as $\theta \uparrow \infty$, which, by $(ii)$ of Assumption \ref{A2}, in turn gives
	$$
	\lim_{\theta \downarrow 0} f(x^{\star}(\theta)) - F^{-1}(\theta) = \infty \quad \text{and} \quad \lim_{\theta \uparrow \infty} f(x^{\star}(\theta)) - F^{-1}(\theta) = - \infty.
	$$
	Then $Q(\widehat{\theta}) > 0$ and $\frac{d}{d\theta}Q(\theta) < 0$ for any $\theta \geq \widehat{\theta}$. 
	Moreover, for any $\theta < \widehat{\theta}$,
	$$Q(\theta) \geq (f(x^{\star}(\theta)) - F^{-1}(\theta))\int_{x^{\star}(\theta)}^{\infty} m'(y) d  y 
	> (f(x^{\star}(\widehat{\theta})) - F^{-1}(\widehat{\theta}))\int_{x^{\star}(\widehat{\theta})}^{\infty} m'(y) d  y = 0,$$
	where the strictly decreasing property of $f \circ x^{\star}$ has been used. Finally, for any $\theta>\theta_o > \widehat{\theta}$ we see that 
	$$
	\frac{d }{d \theta}Q(\theta) 
	< - \frac{1}{F'(F^{-1}(\theta))}\int_{x^{\star}(\theta_o)}^{\infty} m'(y) d  y,
	$$
	where we have used that $x^{\star}(\cdot)$ is decreasing. Hence,
	$$
	Q(\theta) - Q(\theta_o) 
	<  -  \bigg( \int^{\theta}_{\theta_o} \frac{1}{F'(F^{-1}(z))} dz \bigg) \int_{x^{\star}(\theta_o)}^{\infty} m'(y) d  y
	= - (F^{-1}(\theta) - F^{-1}(\theta_0) ) \int_{x^{\star}(\theta_o)}^{\infty} m'(y) d  y,
	$$
	and, taking limits as $\theta\uparrow \infty$ in the latter, and using that $F^{-1}(\theta) \to \infty$, we obtain $Q(\theta) \rightarrow -\infty$. 
	
	All the previous properties of $Q$ imply that there exists a unique $\theta^{\star} > \widehat{\theta}$ solving the consistency equation \eqref{eq:CE}. 
	Therefore, stressing now the dependency of the involved quantities with respect to $r$, and setting $(x_r^\star, \theta_r^\star):=(x^\star(\theta^\star), \theta^\star)$ and $\nu^{\star,r}:=\nu^{\star}(\theta_r^{\star})$, we conclude that  $(\nu^{\star,r}, \theta_r^{\star})$ is the unique equilibrium of the discounted stationary MFG, and that it is characterized by the couple $(x_r^\star, \theta_r^\star)$ solving the system of equations \eqref{eq system discounted}. 
	This completes the proof of the theorem.
	
	
	\subsection{Proof of Theorem \ref{thm ergMFG-ex}}
	\label{sec:proof-Thm-ergMFG}
	
	We divide the proof into two steps. 
	\vspace{0.25cm}
	
	\emph{Step 1.} We fix $\theta>0$ and we solve the control problem with ergodic profit \eqref{eq control proble ergodic}.
	To this aim, define $\mathcal{T}$ as the set of $\mathbb{F}$-stopping times, and, recalling that $b'(x) < -2c$ by $(i)$ of Assumption \ref{assumption ergodic}, consider the auxiliary optimal stopping problem 
	\begin{equation}\label{eq optimal stopping problem}
		u(x,\theta):= \inf_{\tau \in \mathcal{T}} \E_x \bigg[ \int_0^\tau e^{\int_0^t b'(\widehat{X}_s) ds } \pi_x(\widehat{X}_t, \theta) dt  + e^{\int_0^\tau b'(\widehat{X}_s) ds } \bigg].
	\end{equation}
	By employing methods as in \cite{alvarez2001} (see in particular Theorem 5 therein), one can prove that the value function $u(\cdot,\theta)$ is $C^1(\mathbb{R}_+)$ with $u_{xx}(\cdot,\theta) \in L^{\infty}_{\text{loc}}(\R_+)$, and that the optimal stopping time is given by $\tau^\star(x,\theta):=\inf \{t\geq 0 \, | \, \widehat{X}_t^x \leq x^\star(\theta) \} $, where $x^\star(\theta)$ uniquely solves 
	\begin{equation}\label{eq equation for the boundary}
		\int_{x^\star (\theta)}^{\infty} \widehat{\phi}_0 (y) \big(\pi_x(y,\theta) + b'(y) \big) \widehat m'(y) d  y =0.
	\end{equation}
	Existence of a unique solution $x^\star (\theta)$ to the equation \eqref{eq equation for the boundary} can be deduced from Assumption \ref{assumption ergodic} as in Step 1 of the proof of Theorem \ref{thm discMFG-ex}. Also, we have $x^\star(\theta) < \widehat{x}_0(\theta)$ (cf.\ Assumption \ref{assumption ergodic}). Moreover, it can be shown that
	\begin{equation}\label{eq var inequality 1}
		\begin{cases}  
			\mathcal{L}_{\widehat{X}} u(x,\theta)  + b'(x)u(x,\theta) + \pi_x(x,\theta)=0, & x>x^\star(\theta),  \\
			u(x,\theta)=1, & x\leq x^\star(\theta),
		\end{cases}
	\end{equation}
	as well as 
	\begin{equation}\label{eq var inequality 2}
		\begin{cases}
			\mathcal{L}_{\widehat{X}} u(x,\theta)  + b'(x)u(x,\theta) + \pi_x(x,\theta)\geq 0, & x < x^\star(\theta),  \\
			u(x,\theta)\leq 1, & x \geq x^\star(\theta).  
		\end{cases}
	\end{equation}
	Next, define the function $U(\cdot,\theta)$ such that $U_x(x,\theta) = u(z,\theta)$.
	By the regularity of $u(\cdot,\theta)$, the function $U(\cdot,\theta)$ is $C^2(\mathbb{R}_+)$ and we observe that $U_x(x,\theta)=u(x,\theta) \leq 1$ for each $x \in \mathbb{R}_+$.
	Furthermore, setting $\Lambda := b(x^\star(\theta)) + \pi (x^\star(\theta), \theta )$, we find
	\begin{align}\label{eq represent var ineq U}
		\frac{\sigma^2(x)}{2} U_{xx} (x, \theta) + b(x) U_x (x, \theta) + \pi(x,\theta) &=\frac{\sigma^2(x)}{2} u_{x} (x, \theta) + b(x) u (x, \theta) + \pi(x,\theta) \\ \notag
		& = \int_{x^\star(\theta)}^x \Big( \frac{\sigma^2(z)}{2} u_{x} (z, \theta) + b(z) u (z, \theta) + \pi(z,\theta) \Big)_z dz \\ \notag
		& \ \ +\frac{\sigma^2(x^\star(\theta))}{2} u_{x} (x^\star(\theta), \theta) + b(x^\star(\theta)) u (x^\star(\theta), \theta) + \pi(x,\theta) \\ \notag
		& = \int_{x^\star(\theta)}^x \big( \mathcal{L}_{\widehat{X}} u (z, \theta) + b'(z) u (z, \theta) + \pi_x(z,\theta) \big) dz + \Lambda, 
	\end{align}
	where we have used \eqref{eq var inequality 1} and \eqref{eq var inequality 2} in the last equality. 
	Now, if $x< x^\star(\theta)$, the integral in the right-hand side of \eqref{eq represent var ineq U} is nonpositive, so that
	\begin{equation*}
		\frac{\sigma^2(x)}{2} U_{xx} (x, \theta) + b(x) U_x (x, \theta) + \pi(x,\theta) \leq \Lambda. 
	\end{equation*}
	On the other hand, if $x>x^\star(\theta)$, from \eqref{eq represent var ineq U}  and \eqref{eq var inequality 1} we deduce that
	\begin{equation*}
		\frac{\sigma^2(x)}{2} U_{xx} (x, \theta) + b(x) U_x (x, \theta) + \pi(x,\theta) = \Lambda. 
	\end{equation*}
	Overall, we have shown that $U(\cdot,\theta)$ is a $C^2(\mathbb{R}_+)$ function satisfying 
	\begin{equation*}
		\mathcal{L} U(x,\theta) + \pi(x,\theta) \leq \Lambda \quad \text{and} \quad U_x(x,\theta) \leq 1.
	\end{equation*}
	Let now $\nu(x^\star(\theta)) \in \mathcal{A}_e$ be the control that keeps the state process above the threshold  $x^\star(\theta)$. Since $U$ is bounded from below as $U_x(x,\theta)\geq0$ on $\R_+$, and because $\nu(x^\star(\theta))$ increases only when $X^{\nu(x^\star(\theta))} \geq x^\star(\theta)$, a verification theorem shows that $\lambda^\star (\theta)=\Lambda=b(x^\star(\theta)) +  \pi \big(x^\star(\theta), \theta \big)$, and that the process $\nu(x^\star(\theta))$ is optimal. 
	\vspace{0.25cm}
	
	\emph{Step 2.} Given $x^\star(\theta)$ as in Step 1, we impose the consistency condition on $\theta$; that is, we look for $\theta^\star$ such that
	$$
	\int_{x^\star(\theta^\star)}^{\infty} (f(y) - F^{-1}(\theta^\star)) m'(y) d  y = 0.
	$$
	As in Step 2 in the proof of Theorem \ref{thm discMFG-ex}, we can show that such a $\theta^\star$ exists and it is in fact unique. 
	Therefore, setting $(x_e^\star, \theta_e^\star):=(x^\star(\theta^\star), \theta^\star)$ and $\nu^{\star, e}:= \nu(x_e^\star)$ we conclude that $(\nu^{\star,e},\theta_e^\star)$ is the unique equilibrium of the ergodic MFG problem. 
	Moreover, such equilibrium  is characterized by the couple $(x_e^\star, \theta_e^\star)$ uniquely solving \eqref{eq system for the ergodic MFG}, and the value at equilibrium is given by $b(x_e^\star) + \pi(x_e^\star, \theta_e^\star)$. This completes the proof of the theorem.  
	
	
	\subsection{Proof of Lemma \ref{lemma phi continuous}}
	\label{proof:lemma phi continuous}
	
	We prove only the continuity of $\widehat{K}$, the continuity of $\widehat{G}$ being obvious.  Fix $(x,\theta,r) \in \R_+^2 \times (-c,1)$, and a sequence $(x^n,\theta^n,r^n)_{n\in \mathbb{N}}$ converging to $(x,\theta,r)$. 
	Without loss of generality, we can assume that $ a:=x/2<x^n<2x$ and that $\theta/2<\theta^n < 2\theta$ for each $n\in \mathbb{N}$. 
	Also, since the functions $\widehat{\phi}_r$ are defined up to a positive multiplicative factor, we can assume that $\widehat{\phi}_r(a)=1$ for each $r\in (-c,1)$. 
	Hence, for $0<a<y$, defining $\tau_{a}^y:= \inf \{ t \geq 0:\, \widehat{X}_t^{y}\leq a \}$, we have (cf.\ Chapter II in \cite{Borodin Salminen}) 
	\begin{equation}
		\label{eq borodin salminen} 
		\widehat{\phi}_r(y)=\frac{\widehat{\phi}_r(y)}{\widehat{\phi}_r(a)}=\E\bigg[\exp\bigg(\int_0^{\tau_{a}^y}  (b'(\widehat{X}_s^y)-r) ds\bigg) \bigg].
	\end{equation}
	Therefore, for each $0<a<y$ and $-c\leq r< \bar{r} < 1$, one has 
	\begin{equation*}
		\label{eq ratio fundamental solutions}
		{\widehat{\phi}_{\bar{r}}(y)}
		=
		\E\bigg[\exp \bigg( \int_0^{\tau_a^y}  (b'(\widehat{X}_s^y)-\bar{r}) ds\bigg) \bigg] 
		\leq 
		\E\bigg[\exp \bigg( \int_0^{\tau_a^y}  (b'(\widehat{X}_s^y)-r) ds\bigg) \bigg]
		=  
		{\widehat{\phi}_r(y)}, 
	\end{equation*}  
	so that
	\begin{equation}
		\label{eq fundament solution equibounded}
		\widehat{\phi}_1(y) \leq \widehat{\phi}_r(y) \leq \widehat{\phi}_{-c}(y), \quad y>a, \quad r \in (-c,1). 
	\end{equation}
	
	We next prove that $\widehat{\phi}_{r^n}(x^n) \to \widehat{\phi}_r(x)$ as $n \to \infty$. In order to do so, set
	$$
	\alpha^n:=\int_0^{\tau_a^{x^n}}(b'(\widehat{X}_s^{x^n}) -r^n) ds \quad \text{and} \quad \alpha:=\int_0^{\tau_a^{x}}(b'(\widehat{X}_s^{x}) -r) ds, 
	$$
	and observe that ${X}_s^{x^n} \to {X}_s^{x}$ $\P\otimes ds$-a.e.\ and that $\tau_a^{x^n} \to \tau_a^{x}$ $\P$-a.s., as $n\to \infty$. Hence, $$
	\mathds{1}_{(0,\tau_a^{x^n})}(s)(b'(\widehat{X}_s^{x^n}) -r^n) \to \mathds{1}_{(0,\tau_a^{x})}(s)(b'(\widehat{X}_s^{x}) -r), \quad \P\otimes ds\text{-a.s., as } n \to \infty.
	$$
	This, thanks to the Liptschitz continuity of $b$, allows to invoke the dominated convergence theorem in order to deduce that
	\begin{equation}\label{eq limit int tau}
		\alpha^n \to \alpha, \ \P\text{-a.s., as }n \to \infty.
	\end{equation}  
	From \eqref{eq limit int tau} and \eqref{eq borodin salminen}, using that $b'(X_s^{x^n}) - r^n<-c$ for each $n\in \mathbb{N}$, we can employ the dominated convergence theorem once more in order to conclude that
	\begin{equation}\label{eq continuity phi}
		\widehat{\phi}_{r^n}(x^n)=  \E[\exp(\alpha^n)] \to \widehat{\phi}_r(x) = \E[\exp(\alpha^n)], \text{ as $n \to \infty$}.
	\end{equation}
	In the same way, we can prove that
	\begin{equation}\label{eq continuity phi2}
		\widehat{\phi}_{r^n}(y) \to \widehat{\phi}_r(y), \text{ for each $y>a$, as $n \to \infty$}.
	\end{equation}
	
	Next, from \eqref{eq continuity phi} and \eqref{eq continuity phi2}, we have
	\begin{equation}\label{eq limit inside K}
		\mathds{1}_{(x^n,\infty)} (y)\frac{ \widehat{\phi}_{r^n}(y)}{\widehat{\phi}_{r^n}(x^n)} \Pi(y,\theta^n,r^n) \to \mathds{1}_{(x,\infty)}(y)\frac{\widehat{\phi}_{r}(y)}{\widehat{\phi}_{r}(x)}  \Pi(y,\theta,r), 
		\text{ for each $y>a$, as $n\to \infty$.}
	\end{equation}
	Moreover, thanks to \eqref{eq fundament solution equibounded}, we have the estimate 
	$$
	\bigg| \mathds{1}_{(x^n,\infty)} (y)\frac{ \widehat{\phi}_{r^n}(y)}{\widehat{\phi}_{r^n}(x^n)} \Pi(y,\theta^n,r^n)\bigg| 
	\leq \mathds{1}_{(a,\infty)}(y)\frac{\widehat{\phi}_{-c}(y)}{\widehat{\phi}_1(2x)} \Pi(y,\theta/2,-c)  \in \mathbb{L}^1(\R), \quad \text{for each $n\in \mathbb{N}$}. $$
	This, together with \eqref{eq limit inside K}, allows to invoke the dominated convergence theorem and obtain that
	$$
	\widehat{K}(x^n, \theta^n;r^n)
	= \int_{x^n}^\infty \frac{ \widehat{\phi}_{r^n}(y)}{\widehat{\phi}_{r^n}(x^n)} \Pi(y,\theta^n,r^n)dy 
	\to \widehat{K}(x,\theta;r) 
	=\int_{x}^\infty \frac{ \widehat{\phi}_{r}(y)}{\widehat{\phi}_{r}(x)} \Pi(y,\theta,r)dy, \quad
	\text{as $n \to \infty$,}
	$$
	thus providing the claimed continuity of $\widehat{K}$.

	\subsection{Proof of Theorem \ref{thm: abelian}}
	\label{proof:thm: abelian}
	
	We divide the proof in two steps. 
	\vspace{0.25cm}
	
	\emph{Step 1.} In this step we prove the first of the two claimed limits. This is done via a suitable application of the implicit function theorem on the function $\Phi$ (cf.\ \eqref{eq:PHI}), that defines the system of equations characterizing the MFG equilibria. 
	
	For convenience of notation, set $(x^{\star}_0, \theta^{\star}_0):=(x^{\star}_e, \theta^{\star}_e)$.
	Thanks to Lemma \ref{lemma phi continuous}, the map $\Phi$ is continuous and, by Theorem \ref{thm ergMFG-ex}, we have $\Phi(x_0^\star, \theta_0^\star;0)=0$. 
	By invoking Theorem 1.1 in \cite{kumagai80}, the function $r\mapsto (x_r^\star, \theta_r^\star)$ is continuous in a neighborhood $(-\delta, \delta)$ of 0 if and only if there exists neighborhoods $(-\varepsilon, \varepsilon)\subset (-c,1)$ and $B \subset \R^2_+$ of $0$ and of $(x_0^\star,\theta_0^\star)$ respectively, such that the map $\Phi(\cdot,\cdot;r):B \to \R^2$ is locally injective for each $r \in (-\varepsilon,\varepsilon)$. 
	Therefore, we only need to prove local injectivity of the map $\Phi(\cdot,\cdot;r)$, and, in order to accomplish that, we will employ  the local inversion theorem.  In particular, by observing that, for each $r\in (-c,1)$, we have $\Phi(\cdot,\cdot;r)\in C^1(\R_+^2)$, it is enough to show that $\det \mathcal{J} \Phi (x_0^\star,\theta_0^\star;0) \ne 0 $ and that $\det \mathcal{J} \Phi $ is continuous in a neighborhood of $(x_0^\star,\theta_0^\star;0)$, where $\det \mathcal{J}\Phi$ denotes the determinant of the Jacobian matrix of $\Phi$ in the variable $(x,\theta)$.  
	
	We begin by computing the partial derivatives of $\widehat{K}$: 
	\begin{align*}
		{\partial_x \widehat{K}}(x, \theta; r) 
		& = - \Pi(x, \theta; r) - \widehat{K}(x, \theta; r) \frac{\widehat{\phi}_r ' (x)}{\widehat{\phi}_r (x)}  \\ \notag
		&  = - \Pi(x, \theta; r) 
		+ \widehat{K}(x, \theta; r) \widehat{S}'(x) \int_{x}^\infty \frac{\widehat{\phi}_r(y)}{\widehat{\phi}_r(x)} (r-b'(y)) \widehat{m}'(y) dy,  
	\end{align*}
	where, in the second equality, we have used \eqref{psiphiproperties3}. In particular, by repeating arguments similar to those in the proof of Lemma \ref{lemma phi continuous}, one can show that $\partial_x \widehat{K}(x, \theta; r)$ is continuous in $(x,\theta,r)$.
	Also, by Theorem \ref{thm ergMFG-ex}, we have $K(x_0^\star, \theta_0^\star; 0)=0$, so that
	\begin{equation}\label{eq dK dx}
		\partial_x \widehat{K}(x_0^\star, \theta_0^\star; 0) =- \Pi(x_0^\star, \theta_0^\star; 0) < 0, 
	\end{equation}
	where the latter inequality follows from the fact that, arguing as in the proof of Theorem \ref{thm discMFG-ex}, one has
	\begin{equation}
		\label{eq:exRem1}
		\pi_x(x_0^\star, \theta_0^\star) + b'(x_0^\star)>0 
		\quad \text{and} \quad
		f(x_0^\star) -F^{-1}(\theta_0^\star) < 0. 
	\end{equation}
	
	Next, thanks to $\pi$-$(iv)$ in Assumption \ref{A2}, we find
	\begin{equation*}
		\partial_{\theta} \widehat{K}(x, \theta;r)= \int_{x}^\infty \frac{\widehat{\phi}_r(y)}{\widehat{\phi}_r(x) } \big(\pi_{x\theta}(y,\theta) - (r-b'(y)\big) \widehat{m}'(y)  dy, 
	\end{equation*} 
	which, through arguments similar to those in the proof of Lemma \ref{lemma phi continuous}, can be shown to be continuous in $(x,\theta,r)$. Moreover, since $\pi_{x \theta} \leq 0$ and $b' - r < 0$, we have
	\begin{equation}\label{eq dK dtheta}
		\partial_{\theta} \widehat{K}(x_0^\star, \theta_0^\star;0)<0.
	\end{equation} 
	Finally, the function $\widehat{G}$ clearly belongs to $C^1(\R_+^2)$. Moreover,  
	\begin{equation}\label{eq dG dtheta}
		\partial_\theta \widehat{G}(x_0^\star, \theta_0^\star) = - \frac{1}{F'(F^{-1}(\theta_0^\star))}\int_{x_0^\star}^\infty m'(y) dy <0,   
	\end{equation}
	and by \eqref{eq:exRem1} we have
	\begin{equation}\label{eq dG dx}
		\partial_x \widehat{G}(x_0^\star, \theta_0^\star) = - (f(x_0^\star) - F^{-1}(\theta_0^\star) ) m'(x_0^\star) >0. 
	\end{equation}
	
	Therefore, by employing \eqref{eq dK dx}, \eqref{eq dK dtheta}, \eqref{eq dG dtheta} and \eqref{eq dG dx}, and using the continuity of $\det \mathcal{J} \Phi$, we find  neighborhoods $(-\varepsilon, \varepsilon)$  and $B$ of $0$ and $(x_0^\star, \theta_0^\star)$ such that
	\begin{equation*}
		\det \mathcal{J} \Phi (x, \theta;r) = \Big[ \partial_{x} \widehat{K} \,
		\partial_{\theta}\widehat{G} - 
		\partial_{\theta} \widehat{K} \,
		\partial_{x} \widehat{G} \Big] (x, \theta;r) >0, \quad \text{for each $(x,\theta) \in B$ and each $r\in (-\varepsilon, \varepsilon)$}.
	\end{equation*}
	By the latter inequality we can then invoke the local inversion theorem in order to deduce that, for each $r\in (-\varepsilon, \varepsilon)$, the function $\Phi(\cdot,\cdot;r):B \to \R^2$ is locally invertible. Therefore, by Theorem 1.1 in \cite{kumagai80}, the map $r\to (x_r^\star, \theta_r^\star)$ is continuous in $(-\varepsilon, \varepsilon)$.
	\vspace{0.25cm}
	
	\emph{Step 2.}
	With regard to Theorem \ref{thm discMFG-ex} and its proof, we have that
	$$
	rV(x, \theta_r^\star;r) = rV(x^{\star}_r, \theta_r^\star; r) + r(x - x^{\star}_r), \quad x \leq x^{\star}_r.
	$$
	Since $V(\cdot,\theta;r) \in C^2(\mathbb{R}_+)$, by using the fact that $(\mathcal{L}_X - r)V(x, \theta_r^\star;r) + \pi(x,\theta^{\star}_r)=0$ for $x \geq x^{\star}_r$, we find that
	$$ 
	rV(x^{\star}_r, \theta_r^\star; r) = b(x^{\star}_r) + \pi(x^{\star}_r,\theta^{\star}_r).
	$$
	Therefore, for $x \in \mathbb{R}_+$ we can write
	\begin{equation}
		\label{eq:rVr-2}
		rV(x,\theta_r^\star ; r) 
		= r \int_{x^{\star}_r}^{x} V_x(z, \theta_r^\star;r) d z + rV(x^{\star}_r,\theta_r^\star; r)  
		= r \int_{x^{\star}_r}^{x} V_x(z, \theta_r^\star;r) d z + b(x_r^\star) + \pi(x^{\star}_r,\theta^{\star}_r).
	\end{equation}
	Moreover, we have $V_x(\cdot, \theta_r^\star;r) \geq 0$ for each $r >0$. Indeed, for ${z} \leq \bar{z}$ and any control $\nu \in \mathcal{A}_d$, by a comparison theorem (see, e.g., Theorem 54 at p.\ 324 in \cite{protter2005}) we have $X_t^{z; \nu} \leq X_t^{\bar{z}; \nu}$ for each $t \geq 0$, $\P$-a.s. This, together with the monotonicity of $\pi(\cdot, \theta_r^\star)$, implies that 
	$$
	V(z,\theta_r^\star;r)= \sup_{\nu\in \mathcal{A}_d} J(z,\nu, \theta_r^\star;r) \leq \sup_{ \nu\in \mathcal{A}_d} J(\bar{z},\nu, \theta_r^\star;r) = V(\bar{z}, \theta_r^\star;r),
	$$ 
	so that $V_x(\cdot, \theta_r^\star ; r) \geq 0$ for each $r>0$. 
	Also, since $V(\cdot, \theta_r^\star;r)$ solves the equation \eqref{eq HJB discounted} in the proof of Theorem \ref{thm discMFG-ex}, we have $V_x(\cdot, \theta_r^\star;r) \leq 1$ for each $r >0$, which allows to conclude that
	$$
	0 \leq V_x(\cdot, \theta_r^\star;r) \leq 1, \quad \text{ for each $r >0$.}
	$$
	The latter, together with the limits proved in Step 1, allows to use the dominated convergence theorem to take limits as $r\downarrow 0$ in \eqref{eq:rVr-2}, and to conclude that
	$$
	\lim_{r \downarrow 0} r V(x, \theta_r^\star ;r) = \lim_{r \downarrow 0} \Big( b(x_r^\star) + \pi(x^{\star}_r,\theta^{\star}_r)\Big) = b(x_e^\star) + \pi(x^{\star}_e,\theta^{\star}_e)= \lambda^\star (\theta_e^\star),
	$$
	where the last equality follows from Theorem \ref{thm ergMFG-ex}. This completes the proof of the theorem.
	
	
	
	\subsection{Proof of Lemma \ref{lemma a priori estimates}}
	\label{proof:lemma a priori estimates}
	
	We prove only the first estimate, the proof of the second being analogous. Let $i \in \{1,\dots,N\}$ be given and fixed. Recall the definition of $x_{b,\sigma}$, set  $\widehat{L} := \max \{ \widehat{B}, x_{b,\sigma} \}$, and let $\nu \in \mathcal{A}_e(\widehat{B})$. Let then $(\tau^i_k, \bar{\tau}^i_k)_{k\geq1}$ be a sequence of stopping times such that $0\leq \tau^i_1 \leq \bar{\tau}^i_1 \leq \tau^i_2 \leq \bar{\tau}^i_2 \leq \dots $, $\P$-a.s., and such that $\{ X^{i,\nu} \geq \widehat{L} \} = \bigcup_{k \geq 1} [\tau^i_k,\bar{\tau}^i_k]$. 
	
	By employing It\^o's rule on the process $\{ |X^{i,\nu}_t|^2 \}_{t\in [\tau^i_k, \bar{\tau}^i_k]}$, we obtain
	\begin{align*}
		|X^{i,\nu}_{t}|^2 &= |X^{i,\nu}_{\tau^i_k}|^2 +  \int_{\tau^i_k}^t \big( 2X^{i,\nu}_s b (X^{i,\nu}_s) + \sigma^2 (X^{i,\nu}_s) \big) ds  + \int_{\tau^i_k}^t 2 X^{i,\nu}_s \sigma (X^{i,\nu}_s)d W_s  \\ 
		& \leq    |X^{i,\nu}_{\tau^i_k}|^2 +\int_{\tau^i_k}^t 2 X^{i,\nu}_s \sigma (X^{i,\nu}_s) d W_s.
	\end{align*}
	Therefore,
	\begin{align}
		\label{eq:ergestimate}
		\int_0^T \E[ |X^{i,\nu}_t|^2 ] dt 
		& = \int_0^T \E \Big[ |X^{i,\nu}_t|^2 \mathds{1}_{ \{ X^{i,\nu}_t \leq \widehat{L}\} } + |X^{i,\nu}_t|^2 \mathds{1}_{ \{ X^{i,\nu}_t \geq \widehat{L} \} } \Big] dt \nonumber \\
		& \leq \widehat{L}^2 T +  \sum_{k\geq 1}\int_0^T \E \bigg[ \mathds{1}_{(\tau^i_k,\bar \tau^i _k ) }(t) \bigg( |X^{i,\nu}_{\tau^i_k}|^2 +  \int_{\tau^i_k}^t 2 X^{i,\nu}_s \sigma (X^{i,\nu}_s)d W_s   \bigg)   \bigg] dt \nonumber \\
		& \leq  (2 \widehat{L}^2 + \E[ |\xi^i|^2 ])T.
	\end{align}
	In the last inequality above we have used that the expectation of the stochastic integral vanishes and that, because $\nu \in \mathcal{A}_e(\widehat{B})$, one has either $X^{i,\nu}_{\tau^i_k}=\xi^i$ if $\tau^i_k=0$ or $X^{i,\nu}_{\tau^i_k}=\widehat{L}$ if $\tau^i_k>0$.
	
	Since the right-hand side of \eqref{eq:ergestimate} does not depend on the choice of $\nu \in \mathcal{A}_e(\widehat{B})$, the claim is then easily obtained.
	
	
	\subsection{Proof of Theorem \ref{thm approximation N-player game}}
	\label{proof:thm approximation N-player game}
	
	We will prove only Claims 1 and 2, as the proof of Claims 3 and 4 follows similar arguments. 
	\smallbreak\noindent  
	\emph{Proof of Claim 1.} 
	For $\nu \in \mathcal{A}_e(\widehat{B})$, set
	$$
	R^N(\nu):= G^i( \nu, \bar{\boldsymbol{\nu}}^{-i, e}) - G(\nu, \theta_e^\star).  
	$$
	By Theorem \ref{thm ergMFG-ex}, the control policy $\bar\nu^{i,e}$ is optimal for the MFG problem with ergodic cost. Hence
	\begin{equation}\label{eq convergence ergodic rest}
		G^i(\bar{\nu}^{i,e} , \bar{\boldsymbol{\nu}}^{-i, e})  \geq G^i({\nu} , \bar{\boldsymbol{\nu}}^{-i, e})  + R^N(\bar{\nu}^{i,e}) - R^N(\nu),  \quad \nu \in \mathcal{A}_e(\widehat{B}).
	\end{equation}
	Therefore, since $\bar{\nu}^{i,e} \in \mathcal{A}_e(\widehat{B})$ by definition of $\widehat{B}$, we only need to show that $|R^N(\nu)| \to 0$ as $N\to \infty$, uniformly for $\nu \in \mathcal{A}_e(\widehat{B})$. 
	In order to do so, we first observe that 
	\begin{equation}\label{eq rest estimate 1}
		| R^N(\nu) | \leq \limsup_{T \to \infty} \left|\frac{1}{T}\E\left[\int_0^T\pi(X^{i,\nu}_{t},\theta_e^{i,N})dt\right]-\frac{1}{T}\E\left[\int_0^T\pi(X^{i,\nu}_{t},\theta^{\star}_e)dt\right]\right|. 
	\end{equation} 
	Next, using that $\theta_e^{i,N} \geq F(f(x_e^\star))$,  we have by Assumption \ref{ass pi-theta additional}
	\begin{align}\label{eq rest estimate 2} 
		\sup_{\nu  \in \mathcal{A}_e(\widehat{B})} 
		& \limsup_{T \to \infty} \E \left[\int_0^T\frac{1}{T}\left|\pi(X^{i,\nu}_{t},\theta_e^{i,N})-\pi(X^{i,\nu}_{t},\theta^{\star}_e)\right|dt\right] \\ \notag
		& \leq
		\sup_{\nu  \in \mathcal{A}_e(\widehat{B})} \limsup_{T \to \infty} \E \left[\int_0^T\frac{1}{T}C(1+|X^{i,\nu}_{t}|)\left|\theta_e^{i,N}-\theta^{\star}_e\right|dt\right]\\ \notag
		& \leq \E \left[\left|\theta_e^{i,N}-\theta^{\star}_e\right|^2\right]^{\frac{1}{2}} \sup_{\nu  \in \mathcal{A}_e(\widehat{B})} \limsup_{T \to \infty}
		\E \left[\left(\frac{1}{T} \int_0^TC(1+|X^{i,\nu}_{t}|) dt\right)^2\right]^{\frac{1}{2}}\\ \notag
		& \leq \bar{C} \, \E \left[\left|\theta_e^{i,N}-\theta^{\star}_e\right|^2\right]^{\frac{1}{2}}, 
	\end{align} 
	for a constant $ \bar{C}< \infty$ (depending on the initial conditions, but not on $\nu \in \mathcal{A}_e$), and where the last inequality follows from Jensen's inequality and from Lemma \ref{lemma a priori estimates}. 
	
	For $z>0$, set $m_z'(x) := m'(x)\mathds{1}_{[z,\infty)}(x) / \int_z^\infty m'(y) d y$. Exploiting the estimates from Lemma \ref{lemma a priori estimates} and using results from the ergodic theory (see, e.g., p.\ 37 in \cite{Borodin Salminen}), we find
	\begin{align*}
		\E \left[\left|\theta_e^{i,N}-\theta^{\star}_e\right|^2\right] 
		&=  \E \bigg[ \bigg| \lim_{t\to \infty} \frac{1}{t} \int_0^t F \Big( \frac{1}{N-1}\sum_{j\ne i} f(X_s^{j,\bar \nu ^{j,e}}) \Big) ds - \theta_e^\star \bigg|^2 \bigg] \\ 
		&= \bigg|  \int_{\mathbb{R}^{N-1}} F \Big( \frac{1}{N-1}\sum_{j\ne i} f(x^j) \Big) \prod_{j\ne i} m'_{x_e^\star}(x^j)d x^j  - \theta_e^\star \bigg|^2.
	\end{align*} 
	Thanks to the assumption of local Lipschitz continuity of $F$ and the estimates from Lemma \ref{lemma a priori estimates}, we can then employ a suitable version of Hewitt and Savage's theorem (see Corollary 5.13 in \cite{cardaliaguet2010}) and obtain
	$$
	\lim_{N \to \infty} \E \left[\left|\theta_e^{i,N}-\theta^{\star}_e\right|^2\right] 
	= \bigg| F \Big( \int_{\mathbb{R}_+} f(z) m'_{x_e^\star}(z)dz\Big) - \theta_e^\star \bigg|^2 =0.
	$$
	The latter, together with \eqref{eq rest estimate 2} and \eqref{eq rest estimate 1}, gives $|R^N(\nu)|\to 0 $ as $N\to \infty$, uniformly over $\nu \in \mathcal{A}_e(\widehat{B})$. 
	Hence, from \eqref{eq convergence ergodic rest} and \eqref{eq sup = sup}, we conclude the proof of Claim 1.
	
	\smallbreak\noindent 
	\emph{Proof of Claim 2.} Following an argument similar to the one adopted in the previous step, we use Theorem \ref{thm: abelian} and the optimality of $\bar\nu^{i,e}$ for the MFG problem with ergodic cost in order to obtain for any $\nu \in \mathcal{A}_d(\widehat{B})$ that
	\begin{align}
		G^i(\bar{\nu}^{i,r},\bar{\boldsymbol{\nu}}^{-i,r}) 
		& \geq 
		G^i(\nu,\bar{\boldsymbol{\nu}}^{-i,r}) 
		+ G(\nu, \theta_e^\star) - G^i(\nu,  \bar{\boldsymbol{\nu}}^{-i,r}) 
		+ G^i(\bar{\nu}^{i,r}, \bar{\boldsymbol{\nu}}^{-i,r}) - G(\bar{\nu}^{i,r}, \theta_r^\star)
		\\ \notag
		& \quad + G(\bar{\nu}^{i,r}, \theta_r^\star) - G(\bar{\nu}^{i,e},  \theta_r^\star) 
		+ G(\bar{\nu}^{i,e}, \theta_r^\star) - G(\bar{\nu}^{i,e},  \theta_e^\star)\\ \notag 
		& = G^i(\nu,\bar{\boldsymbol{\nu}}^{-i,r}) + G(\bar{\nu}^{i,r}, \theta_r^\star) - G(\bar{\nu}^{i,e},  \theta_r^\star)  + \varepsilon_{N,r},
	\end{align}
	with $\varepsilon_{N,r}$ vanishing as $N \to \infty$ and $r \to 0$.
	Hence, it only remains to show that $G(\bar{\nu}^{i,r}, \theta_r^\star) - G(\bar{\nu}^{i,e},  \theta_r^\star) \to 0$ as $r\to 0$.  
	For $z>0$, set again $m_z'(x) := m'(x) \mathds{1}_{[z,\infty)}(x) / \int_z^\infty m'(y) d y$.
	By the ergodic theory (see, e.g., p.\ 37 in \cite{Borodin Salminen}) and estimates from Lemma \ref{lemma a priori estimates} , for $q\in \{d,e\}$, we have
	\begin{equation} \label{eq ergodic thm pi}
		\lim_{T\to\infty}\frac{1}{T}\E\left[\int_0^T\pi(X^{i,\bar{\nu}^{i,q}}_{t},\theta^\star_r )dt\right]
		=\int_{x_{q}^\star}^\infty \pi\left(x,\theta^\star_r\right) m'_{x_q^\star}(x)dx.
	\end{equation}
	Also, as in the proof of Lemma \ref{lemma a priori estimates}, one can prove that $\E [ X_T^{i,\bar{\nu}^{i,q}} ] \leq 2 \widehat{L} + \E [\xi^i]$, and we obtain
	\begin{equation}\label{eq ergodic control}
		\lim_{T\to \infty} \frac{1}{T} \E [\bar{\nu}^{i,q}_T] = \lim_{T\to \infty} \frac{1}{T} \E \bigg[ X_T^{i,\bar{\nu}^{i,q}}- \xi^i - \int_0^T b(X_s^{i,\bar{\nu}^{i,q}})ds \bigg] =  - \int_{x_{q}^\star}^\infty b(x) m'_{x_q^\star}(x)dx.
	\end{equation}
	
	
	By the convergence in Theorem \ref{thm: abelian}, using \eqref{eq ergodic thm pi} and \eqref{eq ergodic control}, we conclude that
	$G(\bar{\nu}^{i,r}, \theta_r^\star) - G(\bar{\nu}^{i,e},  \theta_r^\star) \to 0$ as $r\to 0$, thus completing the proof of Claim 2.

	\medskip 
	

%
%



\textbf{Acknowledgments} {Jodi Dianetti and Giorgio Ferrari acknowledge financial support from the German Research Foundation (DFG) through the CRC 1283. The authors also thank Xin Guo for many insightful discussions.}


\begin{thebibliography}{99}


\bibitem[{Achdou et~al.(2014)}]{Achdouetal} Achdou Y, Buera FJ, Lasry J-M, Lions J-P, Moll B (2014) Partial differential equations models in macroeconomics. {\it Phil.\ Trans.\ R.\ Soc.\ A} 372(2028).

\bibitem[{Adlakha et~al.(2015)}]{Adlakhaetal} Adlakha S, Johari R, Weintraub GY\ (2015) {Equilibria of dynamic games with many players: Existence, approximation, and market structure}. {\it J.\ Econ.\ Theory} 156:269--316.

\bibitem[{A\"id et~al.(2020)}]{Tankov2} A\"id R, Dumitrescu R, Tankov P (2020) {The entry and exit game in the electricity markets: A mean field game approach}. Preprint, arXiv: 2004.14057.

\bibitem[{Alvarez(2001)}]{alvarez2001} Alvarez LHR (2001) Reward functionals, salvage values, and optimal stopping. {\it Math.\ Meth.\ Oper.\ Res.}\ 54(2):315-337.

\bibitem[{Alvarez(2015)}]{Alvarez} Alvarez LHR  (2015) {Singular stochastic control, linear diffusions, and optimal stopping: A class of solvable problems}. {\it SIAM J.\ Control Optim.}\ 39(6):1697--1710.

\bibitem[{Alvarez and Matom\"aki(2015)}]{AlvarezMat} Alvarez LHR,  Matom\"aki P (2015) {Expected supremum representation and optimal stopping}. Preprint, arXiv: 1505.01660.

\bibitem[Alvarez(2018)]{Alvarez2018} Alvarez LHR {A Class of solvable stationary singular stochastic control problems}. Preprint, arXiv: 1803.03464.


\bibitem[{Alvarez and Hening(2019)}]{AlvarezHening} Alvarez LHR, Hening A (2019) Optimal sustainable harvesting of populations in random environments. {\it Stoch.\ Proc.\ Appl.}\ forthcoming.

\bibitem[{Baldursson and Karatzas(1997)}]{BalKar} Baldursson FM, Karatzas I (1997) {Irreversible investment and industry equilibrium}. {\it Fin.\  Stoch.}\ 1:69--89.

\bibitem[{Bardi(2016)}]{Bardi} Bardi M (2016) {Nonlinear elliptic systems and mean-field games}. {\it Nonlinear Differ.\ Equ.\ Appl.}\ 23:44.

\bibitem[{Basei et~al.(2020)}]{Baseietal} Basei M, Cao H, Guo X (2020) {Nonzero-sum stochastic games and mean field games with impulse controls}. {\it Math.\ Oper.\ Res.}\ forthcoming.

\bibitem[{Bensoussan et~al.(2010)}]{Bensoussanetal} Bensoussan A, Liu J, Yuan J (2010) {Singular control and impulse control: A common approach}. {\it Discrete Cont.\ Dyn.-B} 13(1):27--57.

\bibitem[{Bensoussan et~al.(2016)}]{Bensoussan-etal-2016} Bensoussan A, Sung KCJ, Yam SCP, Yung SP (2016). Linear-quadratic mean field games. {\it J.\ Optim.\ Theory Appl.}\ 169:496--529.

\bibitem[{Bertola(1998)}]{Bertola} Bertola G (1998). Irreversible investment. {\it Res. Econ.}\ 52(1):3--37.

\bibitem[{Bertucci(2018)}]{Bertucci} Bertucci C (2018) Optimal stopping in mean field games, an obstacle problem approach. {\it J.\ Math.\ Pures Appl.}\ 120:165--194.


\bibitem[{Borodin and Salminen(2015)}]{Borodin Salminen} Borodin AN, Salminen, P (2015) {\it Handbook of Brownian motion - Facts and formulae} (Springer).

\bibitem[{Bouveret et~al.(2020)}]{Tankov1} Bouveret G, Dumitrescu R, Tankov P (2020) {Mean field games of optimal stopping: a relaxed solution approach}. {\it SIAM J.\ Control Optim.}\ 58(4):1795--1821.

\bibitem[{Campi et~al.(2020)}]{Campietal} Campi L, De Angelis T, Ghio M, Livieri G (2020) {Mean field games of finite-fuel capacity expansion with singular controls}. Preprint, arXiv: 2006.02074.

\bibitem[{Cao and Guo(2020)}]{CaoGuo} Cao H, Guo X (2020) {MFGs for partially reversible investment}. {\it Stoch.\ Proc.\ Appl.}\ forthcoming.

\bibitem[{Cardaliaguet(2010)}]{cardaliaguet2010} Cardaliaguet, P (2010). {Notes on mean field games}. {Technical report}.


\bibitem[{Carmona and Delarue(2018)}]{CarmonaDelarue18} Carmona R, Delarue F (2018) {\it Probabilistic theory of mean field games with applications} (Springer).

\bibitem[{Christensen et~al.(2020)}]{Christensen} Christensen S, Neumann BA, Sohr T (2020) {Competition versus cooperation: A class of solvable mean field impulse control problems}. Preprint, arXiv: 2010.06452.

\bibitem[{Cirant(2016)}]{Cirant} Cirant M (2016) {Stationary focusing mean-field games}. {\it Commun.\ Partial\ Differ.\ Equ.}\ 41(8):1324--1346.

\bibitem[{De Angelis and Ferrari(2018)}]{DeAFe} De Angelis T, Ferrari G (2018). Stochastic nonzero-sum games: A new connection between singular control and optimal stopping. {\it  Adv.\ Appl.\ Probab.}\ 50(2):347--372.

\bibitem[{Dianetti and Ferrari(2020)}]{DianettiFerrari} Dianetti J, Ferrari G (2020). Nonzero-sum submodular monotone follower games: Existence and approximation of Nash equilibria. {\it SIAM J.\ Control Optim.}\ 58(3):1257--1288.

\bibitem[{Dixit and Stiglitz(1977)}]{DixitStigl} Dixit AK, Stiglitz JE (1977). Monopolistic competition and optimum product diversity. {\it Am.\ Econ.\ Rev.}\ 67(3):297--308.

\bibitem[{Dixit and Pindyck(1994)}]{DixitPindyck} Dixit AK, Pindyck RS (1994). {\it Investment under uncertainty}. Princeton University Press. 

\bibitem[{Ferrari et al.(2017)}]{FeStRi} Ferrari G, Steg J-H, Riedel F (2017). Continuous-time public good contribution under uncertainty: A stochastic control approach. {\it Appl.\ Math.\ Optim.}\ 75(3):429--470.


\bibitem[{Fudenberg and Tirole(1991)}]{Fudenberg-Tirole} Fudenberg D, Tirole J (1991) {\it Game theory} (MIT Press).

\bibitem[{Gabaix(2009)}]{Gabaix} Gabaix X (2009) Power laws in Economics and Finance {\it Annu.\ Rev.\ Econ.}\ 1(1):255--293.

\bibitem[{Gu\'eant et~al.(2011)}]{Gueant-etal2011} {Gu\'eant O, Lasry J-M, Lions P-L} (2011) {\it Mean field games and applications}, Paris-Princeton Lectures on Mathematical Finance 2010, volume 2003 of Lecture Notes in Mathematics, pp.\ 205--266 (Springer).

\bibitem[{Guo and Xu(2019)}]{GuoXu} Guo X, Xu R (2019) {Stochastic games for fuel follower problem: $N$ versus mean field game}. {\it SIAM J.\ Control Optim.}\ 57(1):659--692.

\bibitem[{Guo et~al.(2020)}]{Guoetal} Guo X, Tang W, Xu, R (2020) A class of stochastic games and moving free boundary problems. Preprint, arXiv: 1809.03459.

\bibitem[{Harrison(2013)}]{Harrison} Harrison MJ (2013) {\it Brownian models of performance and control} (Cambridge University Press).

\bibitem[{Hopenhayn(1992)}]{Hopenhayn} Hopenhayn HA (1992) {Entry, exit and firm dynamics in long run equilibrium}. {\it Econometrica} 60(5):1127--1150.

\bibitem[{Horst and Fu(2017)}]{HorstFu} Horst U, Fu G (2017) {Mean field games with singular controls}. {\it SIAM J.\ Control Optim.}\ 55(6):3833--3868.

\bibitem[{Huang et~al.(2006)}]{HuangMalhameCaines06} Huang M, Malham{\'e} RP, Caines PE (2006) {Large population stochastic dynamic games: Closed-loop McKean-Vlasov systems and the Nash certainty equivalence principle} {\it Commun.\ Inf.\ Syst.}\ {6(3)}:221-252.

\bibitem[{Jack and Zervos(2006)}]{JackZervos} Jack A, Zervos M (2006) {A Singular control problem with an expected and a pathwise ergodic performance criterion}. {\it J.\ Appl.\ Math.\ Stoch.\ Anal.}\ Article ID 82538:1--19.

\bibitem[{Jack et al.(2008)}]{JackJonhnsonZervos} Jack A, Jonhnson TC, Zervos M (2008) {A singular control problem with application to the goodwill problem}. {\it Stoch.\ Proc.\ Appl.}\ 118:2098--2124.

\bibitem[{Jovanovic and Rosenthal(1988)}]{JovanovicRos} Jovanovic B, Rosenthal R W (1988) {Anonymous sequential games} {\it J.\ Math.\ Econ.}\ {17}:77--87. 

\bibitem[{Karatzas(1983)}]{Kar83} Karatzas I (1983) {A class of singular stochastic control problems}. {\it Adv.\ Appl.\ Prob.}\ 15:225--254.

\bibitem[{Karatzas and Shreve(1998)}]{KS} Karatzas I, Shreve SE (1998) {\it Brownian motion and stochastic calculus} (Springer).

\bibitem[{Kumagai(1980)}]{kumagai80} Kumagai S (1980) {An implicit function theorem: Comment.} {\it J.\ Optim.\ Theory Appl.}\ 31(2):285--288.

\bibitem[{Kwon(2020)}]{Kwon} Kwon HD (2020) {Game of variable contributions to the common good under uncertainty} {\it Oper.\ Res.}\ forthcoming.

\bibitem[{Kwon and Zhang(2015)}]{KwonZhang} Kwon HD, Zhang H (2015). Game of singular stochastic control and strategic exit. {\it Math.\ Oper.\ Res.}\ 40(4):869--887.

\bibitem[{Lasry and Lions(2007)}]{LasryLions07} Lasry J-M, Lions P-L (2007) Mean field games. {\it Japanese J.\ Math.}\ 2(1):229--260.

\bibitem[{Lon and Zervos(2011)}]{LZ11} Lon PC, Zervos M (2011) {A model for optimally advertising and launching a product}. {\it Math.\ Oper.\ Res.}\ 36:363--376.

\bibitem[{Luttmer(2007)}]{Luttmer} Luttmer EGJ (2007) {Selection, growth, and the size distribution of firms}. {\it Q.\ J.\ Econ.} 122(3):1103--1144.

\bibitem[{Miao(2005)}]{Miao} Miao J (2005) Optimal capital structure and industry dynamics. {\it J.\ Finance} 60(6):2621--2659.


\bibitem[{Pilipenko(2014)}]{Pilipenko2014} Pilipenko A (2014) An introduction to stochastic differential equations with reflection. Universit\"atsverlag Potsdam.

\bibitem[{Protter(2005)}]{protter2005} Protter, PE (2005) {\it Stochastic integration and differential equations} (Springer).

\bibitem[{Weintraub et~al.(2008)}]{Weintraubetal2008} Weintraub GY, Benkard CL, Van Roy B (2008) Markov perfect industry dynamics with many firms. {\it Econometrica} {76(6)}:1375--1411.

\bibitem[{Weintraub et~al.(2011)}]{Weintraubetal2011} Weintraub GY, Benkard CL, Van Roy B (2011). {Industry dynamics: Foundations for models with an infinite number of firms}. {\it J.\ Econ.\ Theory} {46}:1965--1994.

\bibitem[{Weerasinghe(2007)}]{Weerasinghe2007} Weerasinghe A (2007). {An Abelian limit approach to a singular ergodic control problem}. {\it SIAM J.\ Control\ Optim.}\ {46(2)}:714--737.




\end{thebibliography}
\end{document}